\documentclass{article}

\usepackage{amssymb}
\usepackage{amsthm}
\usepackage{amscd}
\usepackage{amsmath}
\theoremstyle{plain}
\newtheorem{theorem}{Theorem}[section]
\newtheorem*{theorem*}{Theorem}

\newtheorem*{maintheorem*}{Main Theorem}
\newtheorem{proposition}[theorem]{Proposition}
\newtheorem{corollary}[theorem]{Corollary}
\newtheorem{lemma}[theorem]{Lemma}
\newtheorem{claim}[theorem]{Claim}
\newtheorem*{conjecture*}{Conjecture}

\theoremstyle{definition}
\newtheorem{definition}[theorem]{Definition}
\newtheorem*{definition*}{Definition}
\newtheorem{example}[theorem]{Example}
\newtheorem*{example*}{Example}
\newtheorem{notation}[theorem]{Notation}
\newtheorem*{notation*}{Notation}
\newtheorem*{notation-conv*}{Notation and convention}
\newtheorem*{convention*}{Convention}
\theoremstyle{remark}

\newcommand{\Z}{{\mathbb Z}}
\newcommand{\C}{{\mathbb C}}
\newcommand{\Q}{{\mathbb Q}}
\newcommand{\R}{{\mathbb R}}
\newcommand{\proj}{{\mathbb P}}
\newcommand{\K}{{\mathbb K}}

\begin{document}
\title{Fourier-Mukai Partners of a $K3$ Surface and the Cusps of its Kahler Moduli} 
\author{Shouhei Ma \\
{\small Graduate School of Mathematical Sciences, University of Tokyo   } \\
{\small 3-8-1\: Komaba\: Meguro-ku, Tokyo\: 153-8914, Japan } \\
{\small E-mail: sma@ms.u-tokyo.ac.jp } }
\date{}
\maketitle

\begin{abstract}
Using lattice theory, 
we establish 
a one-to-one correspondence 
between 
the set of Fourier-Mukai partners of a projective $K3$ surface 
and 
the set of $0$-dimensional standard cusps of its Kahler moduli. 
We also study 
the relation 
between 
twisted Fourier-Mukai partners 
and general $0$-dimensional cusps, 
and the relation 
between 
Fourier-Mukai partners with elliptic fibrations 
and 
certain $1$-dimensional cusps.
\end{abstract}

\section{Introduction}
The excellent work of Mukai (\cite{Mu1}, \cite{Mu2}) and Orlov (\cite{Or}) 
opened a way to study derived equivalence for $K3$ surfaces via their periods. 
One of their conclusions is that 
two projective $K3$ surfaces are derived equivalent if and only if 
there exists an isometry between their Mukai lattices preserving the periods.
This theorem can be viewed as a generalization of 
the global Torelli theorem for $K3$ surfaces. 
Let ${\rm FM\/}(S)$ be the set of isomorphism classes of Fourier-Mukai (FM) partners of 
a projective $K3$ surface $S$, i.e., 
$K3$ surfaces derived equivalent to $S$.
As an application of Mukai-Orlov's theorem, 
Hosono-Lian-Oguiso-Yau (\cite{H-L-O-Y}) gave a counting formula for 
$\# {\rm FM\/}(S)$.

In this paper, 
we construct 
a bijection between
${\rm FM\/}(S)$ and 
the set of embeddings of the hyperbolic plane $U$ 
into the lattice $\widetilde{NS}(S):=NS(S)\oplus U$ 
up to the action of a certain finite-index subgroup 
$\Gamma _{S}\subset O(\widetilde{NS}(S))$. 
Then, 
by considering the isotropic vector $(0, 1)\in U$, 
we obtain 
$0$-dimensional standard cusps of the modular variety 
$\Gamma _{S}^{+}\backslash \Omega _{\widetilde{NS}(S)}^{+}$ 
 associated with 
$\Gamma _{S}^{+}$ and $\widetilde{NS}(S)$ (\cite{B-B}, \cite{PS}).
Let 
$\Gamma _{S}^{+}\backslash I^{1}(\widetilde{NS}(S))$ 
be the set of $0$-dimensional standard cusps of 
$\Gamma _{S}^{+}\backslash \Omega _{\widetilde{NS}(S)}^{+}$.
We shall prove the following.

\begin{theorem}[Theorem \ref{main1}]
There exists a bijective map 
\[
\mu _{0} : {\rm FM\/}(S)\rightarrow \Gamma _{S}^{+}\backslash I^{1}(\widetilde{NS}(S)).
\]
In particular,
\[
\# {\rm FM\/}(S)=\# 
\{ \text{ 0-dimensional standard cusps of } \Gamma _{S}^{+}\backslash \Omega _{\widetilde{NS}(S)}^{+} \} .
\]
\end{theorem}

Of course, 
not every $0$-dimensional cusp of 
$\Gamma _{S}^{+}\backslash \Omega _{\widetilde{NS}(S)}^{+} $
is standard.
It turns out that 
the existence of non-standard cusps 
corresponds to 
the existence of {\it twisted\/} FM-partners of $S$.
Let 
${\rm FM\/}^{d}(S)$
be the set of isomorphism classes of twisted FM-partners $(S', \alpha ')$ of $S$ with ${\rm ord\/}(\alpha ')=d$, 
and let
$\Gamma _{S}^{+}\backslash I^{d}(\widetilde{NS}(S))$ 
be the set of $0$-dimensional cusps $[l]$ of 
$\Gamma _{S}^{+}\backslash \Omega _{\widetilde{NS}(S)}^{+}$ 
with ${\rm div\/}(l)=d$. 
Introducing certain quotient sets
$\mathcal{FM}^{d}(S)$,  
$r(\Gamma _{S}^{+})\backslash I^{d}(A_{\widetilde{NS}(S)})$ 
of 
${\rm FM\/}^{d}(S)$,  
$\Gamma _{S}^{+}\backslash I^{d}(\widetilde{NS}(S))$
respectively 
(see Sect.3.2 for the definitions), 
and using 
Huybrechts-Stellari's solution of C\u ald\u araru conjecture (\cite{H-S2}),
we shall prove the following.

\begin{theorem}[Theorem \ref{main2}]
There exist a map 
\[
\nu _{0} :  \Gamma _{S}^{+}\backslash I^{d}(\widetilde{NS}(S)) 
\longrightarrow {\rm FM\/}^{d}(S)
\]
and a bijective map 
\[
\xi _{0} : r(\Gamma _{S}^{+})\backslash I^{d}(A_{\widetilde{NS}(S)})
\stackrel{\simeq}{\longrightarrow } 
\mathcal{FM}^{d}(S)
\]
which fit in the following commutative diagram:
$$\CD
\Gamma _{S}^{+}\backslash I^{d}(\widetilde{NS}(S)) @>\nu _{0}>> {\rm FM\/}^{d}(S) \\
@VpVV   @VV\pi V \\
r(\Gamma _{S}^{+})\backslash I^{d}(A_{\widetilde{NS}(S)}) @>>\xi _{0}> \mathcal{FM}^{d}(S).
\endCD $$
\end{theorem}

Via Theorem 1.1 and 1.2, 
we can obtain informations 
about the abstract set $ {\rm FM\/}^{d}(S)$ 
by studying 
the $0$-dimensional cusps of the modular variety 
$\Gamma _{S}^{+}\backslash \Omega _{\widetilde{NS}(S)}^{+} $.

Besides $0$-dimensional cusps, 
$\Gamma _{S}^{+}\backslash \Omega _{\widetilde{NS}(S)}^{+} $
has also $1$-dimensional cusps.
Let ${\rm FM\/}_{ell}(S)$ be the set of isomorphism classes of 
pairs $(S', L')$, where $S'\in {\rm FM\/}(S)$ and 
$L'\in Pic(S')$ is 
the line bundle associated to a smooth elliptic curve on $S'$.
In the same manner as 
in the case of $0$-dimensional standard cusps, 
we relate ${\rm FM\/}_{ell}(S)$ to
the set $\Gamma _{S}^{+}\backslash I_{2}^{st}(\widetilde{NS}(S))$ 
of those $1$-dimensional cusps of 
$\Gamma _{S}^{+}\backslash \Omega _{\widetilde{NS}(S)}^{+}$ 
whose closure contains a $0$-dimensional standard cusp.
The map $\mu _{1}$ 
from ${\rm FM\/}_{ell}(S)$ to
the set of such $1$-dimensional cusps 
is surjective but not injective in general.
For example, 
for generic $K3$ surface $S$ with $NS(S)=U(r)$ $(r>2)$, 
we will calculate $\mu _{1}$ explicitly
to conclude that 
$\mu _{1}$ is far from being injective (Example \ref{U(r)}).
On the other hand, 
if ${\rm det\/}NS(S)$ is square-free, 
then $\mu _{1}$ is bijective (Corollary \ref{inj if sq-free}).
It turns out that 
${\rm FM\/}_{ell}(S)$ carries informations about 
the compactifications of $1$-dimensional cusps which belong to  
$\Gamma _{S}^{+}\backslash I_{2}^{st}(\widetilde{NS}(S))$.
We shall observe that 
if two elliptic $K3$ surfaces $(S,L)$ and $(S',L')$ 
give the same $1$-dimensional cusp, 
there exists a coherent sheaf on $S\times S'$ 
via which $L$ gives rise to $L'$. 
(For precise statement see Proposition \ref{FM transform} and the remark after it.)

According to Bridgeland's results (\cite{Br1}, \cite{Br2}), 
$\Gamma _{S}^{+}\backslash \Omega _{\widetilde{NS}(S)}^{+} $
is biholomorphic  to a natural quotient of 
the space of certain stability conditions on $D^{b}(S)$, 
after removing some divisors on 
$\Omega _{\widetilde{NS}(S)}^{+} $.
That is, we have a canonical isomorphism (Proposition \ref{Bridgeland})
\[ 
Aut^{\dag }(D^{b}(S))\backslash 
Stab^{\dag }(S)/\widetilde{GL}_{2}^{+}({\R}) 
\simeq
\Gamma _{S}^{+} \backslash 
\Bigl( \Omega _{\widetilde{NS}(S)}^{+}-\bigcup_{\delta \in \Delta (S)}\delta ^{\perp} \Bigr) .
\] 
Hence 
Theorem 1.1 asserts that 
${\rm FM\/}(S)$ is identified with 
the $0$-dimensional standard cusps of 
the Baily-Borel compactification of 
$Aut^{\dag }(D^{b}(S))\backslash 
Stab^{\dag }(S)/\widetilde{GL}_{2}^{+}({\R})$.
However, 
the lattice theoretic approach 
employed in this paper 
does not explain any intrinsic reason 
for the correspondence between 
(twisted) FM-partners 
and $0$-dimensional cusps.
It may be interesting to 
understand the correspondence  
in more derived-categorical way.

This paper is organized as follows.
In Sect.2.1, 
we recall some facts about even lattices and their discriminant forms, 
following Nikulin (\cite{Ni}).
Sect.2.2 is devoted to 
the study of isotropic elements of an even lattice. 
In Sect.2.3, 
we recall 
the Baily-Borel compactification of 
an arithmetic quotient of 
a type IV symmetric domain, 
following \cite{B-B}, \cite{PS}, \cite{Sc}.
In Sect.3.1, 
we prove Theorem 1.1.
In Sect.3.2, 
we prove Theorem 1.2. 
In Sect.4.1, 
we define ${\rm FM\/}_{ell}(S)$ and relate it to 
primitive embeddings of $U\oplus {\Z}l$ into $\widetilde{NS}(S)$.
In Sect.4.2, we study the relation between ${\rm FM\/}_{ell}(S)$ and certain $1$-dimensional cusps of 
$\Gamma _{S}^{+}\backslash  \Omega _{\widetilde{NS}(S)}^{+}$.
In Sect.5, 
we observe the action of 
$\Gamma _{S}^{+}$ on 
$\Omega _{\widetilde{NS}(S)}^{+}$,
and 
we also describe $\Gamma _{S}^{+}\backslash  \Omega _{\widetilde{NS}(S)}^{+}$ 
in terms of the 
space of stability conditions.

\begin{notation}
By an $even$ $lattice$, we mean a free ${\mathbb Z}$-module $L$ of finite rank 
equipped with 
a non-degenerate symmetric bilinear form $( , ):L\times L\to {\Z}$ 
satisfying $(l, l)\in 2{\Z}$ for all $l\in L$. 
Denote by ${\rm rk\/}(L)$ and ${\rm sign\/}(L)$ 
the rank and the signature of $L$, respectively.
For a lattice $L$ and a field ${\K}$,  
$L_{{\K}}$ denotes the ${\K}$-vector space $L\otimes {\K} $.
For two lattices $L$ and $M$, 
$L\oplus M$ is the lattice defined as 
the orthogonal direct sum of $L$ and $M$,
while $L+M$ denotes the direct sum of 
the ${\Z}$-modules underlying $L$ and $M$.
The projection from $L\oplus M$ to $L$ is denoted by 
$pr_{L}:L\oplus M \rightarrow L$.
The group of isometries of $L$ is denoted by $O(L)$. 
For an element $l \in L$, 
we define the positive integer 
${\rm div\/}(l)$ 
to be the generator of the ideal $(x, L) \subseteq {\mathbb Z}$. 
A sublattice $M\subset L$ is called $primitive$ 
if $L/M$ is a free ${\Z}$-module. 
For (possibly degenerate) lattices $L$ and $M$,  
we denote by ${\rm Emb\/}(L, M)$ the set of primitive embeddings of $L$ into $M$. 
A sublattice $M \subseteq L$ is called $isotropic$ 
if $(x, y)= 0$ 
for all $x, y\in M$. 
We denote by $I_{r}(L)$  
the set of primitive isotropic sublattices of $L$ of ${\rm rk\/}=r$. 
A non-zero element $l\in L$ is called isotropic 
(resp. primitive) 
if ${\Z}l$ is isotropic 
(resp. primitive).
We denote  
$I^{d}(L):=\{ l\in L | \: l \: \text{is primitive }, \: (l, l)=0, \: {\rm div\/}(l)=d \: \} $.

By a $K3$ surface, 
we mean a projective $K3$ surface over ${\C}$.
For a $K3$ surface $S$, 
we denote by $NS(S)$ (resp. $T(S)$) 
the Neron-Severi (resp. transcendental) lattice of $S$. 
Set 
$\widetilde{NS}(S):=H^{0}(S, {\Z}) + NS(S) + H^{4}(S, {\Z})$, 
$\Lambda _{K3}:=U^{3}\oplus E_{8}^{2}$, 
and 
$\widetilde{\Lambda }_{K3}:=U^{4}\oplus E_{8}^{2}$. 
\end{notation}

\section{Preliminaries from lattice theory}
\subsection{Even lattice and its discriminant form}

For an even lattice $L$, 
we can associate a finite Abelian group 
$A_{L}:=L^{\vee}/L$ 
equipped with a quadratic form 
$q_{L}:A_{L}\rightarrow {\mathbb Q} / 2\mathbb{Z}$, 
$q_{L}(x+L) = (x, x) + 2{\mathbb Z}$.
The corresponding bilinear form on $A_{L}$ is 
$b_{L}:A_{L}\times A_{L}\rightarrow {\mathbb Q}/{\mathbb Z}$, 
$b_{L}(x+L, y+L)=(x, y)+{\mathbb Z}$.
We call $(A_{L},q_{L})$ or $(A_{L},b_{L})$ the {\it discriminant form\/} of $L$. 
There is a natural homomorphism 
$r_{L}:O(L)\rightarrow O(A_{L}, q_{L})$, 
whose kernel is denoted by $O(L)_{0}\subset O(L)$. 
We often write just 
$A_{L}$ or $(A_{L}, q)$ 
instead of $(A_{L}, q_{L})$, 
and write $r$ instead of $r_{L}$.  
The following well-known facts due to Nikulin  will be used frequently in this paper :

\begin{proposition}[\cite{Ni}]\label{Nikulin1}
Let $M$ be a primitive sublattice of an even unimodular lattice $L$ 
with the orthogonal complement $M^{\perp}$. 
Then

{\rm (1)\/} There is an isometry 
$\lambda _{L}: 
(A_{M},q_{M})\stackrel{\simeq}{\rightarrow} (A_{M^{\perp}},-q_{M^{\perp}})$.

{\rm (2)\/} For isometries 
$\gamma _{M} \in O(M)$ and $\gamma _{M^{\perp}} \in O(M^{\perp})$, 
there is an isometry $\gamma _{L}\in O(L)$ such that 
$\gamma _{L}|_{M\oplus M^{\perp}}=\gamma _{M} \oplus \gamma _{M^{\perp}}$ 
if and only if 
$\gamma _{M} $ and $\gamma _{M^{\perp}} $ are compatible on the discriminant group, i.e, 
$\lambda _{L}\circ r_{M}(\gamma _{M} )=r_{M^{\perp}}(\gamma _{M^{\perp}} )\circ \lambda _{L}$.
\end{proposition}

Two even lattices $L$ and $M$ are said to be $isogenus$ if 
$L\otimes {\Z}_{p} \simeq M\otimes {\Z}_{p}$ for every prime number $p$  
and ${\rm sign\/}(L)={\rm sign\/}(M)$ .
By \cite{Ni}, this is equivalent to the condition that 
$(A_{L},q_{L})\simeq (A_{M},q_{M})$ and ${\rm sign\/}(L)={\rm sign\/}(M)$ . 
The set of isometry classes of lattices isogenus to $L$ is denoted by $\mathcal{G} (L)$.

\begin{proposition}[\cite{Ni}]\label{Nikulin2}
Let $l(A_{L})$ be the minimal number of the generators of the Abelian group $A_{L}$.
If $L$ is indefinite and \, ${\rm rk\/}(L)\geq l(A_{L})+2$, \, 
then $\mathcal{G} (L)= \{ L\} $ 
and the homomorphism $r_{L}$ is surjective.
\end{proposition}

By Proposition \ref{Nikulin2}, $L$ and $M$ are isogenus if and only if $L\oplus U\simeq M\oplus U$,  
where $U$ is the even indefinite unimodular lattice ${\Z}e+{\Z}f,(e,f)=1,(e,e)=(f,f)=0$.
Note that if an even unimodular lattice $L$ is embedded in another even lattice $M$, 
$L$ must be an orthogonal summand of $M$.

\begin{lemma}\label{proj}
Let $L$ and $M$ be even lattices. 
If there is an isometry $\varphi :L\oplus U\simeq M\oplus U$ 
with $\varphi (f)=f$, then the composition  
$pr_{M}\circ (\varphi |_{L}) : L\hookrightarrow \varphi (f)^\perp =M\oplus {\Z}f
\twoheadrightarrow M$
is an isometry.
\end{lemma}

\begin{proof}
Since 
$L\oplus {\Z}f=f^{\perp}\simeq \varphi (f)^{\perp}=M\oplus {\Z}f$, 
we can write $\varphi (l)=m+\alpha f$ for each $l\in L$, 
where 
$m\in M$ and $\alpha \in {\Z}$.
Then $(l,l)=(\varphi (l),\varphi (l))=(m,m)$. 
Thus the correspondence $\varphi (l)\mapsto m$ gives 
an embedding of lattice $L\hookrightarrow M$.
By the inclusions $L\subseteq M\subseteq M^{\vee }\subseteq L^{\vee }$, 
we have 
$|A_{M}|=|M/M^{\vee}|\leq |L/L^{\vee}|=|A_{L}|$.
Here the equality holds if and only if $L\simeq M$ by 
$pr_{M}\circ (\varphi |_{L})$.
Since $A_{L}\simeq A_{L\oplus U}\simeq A_{M\oplus U}\simeq A_{M}$,
we get $|A_{M}|=|A_{L}|$, 
so that $L\simeq M$.
\end{proof}

\subsection{Primitive isotropic vectors}

Let $L$ be an even lattice possessing a primitive isotropic vector $l$. 
We shall study some properties of $L$ related to $l$.

\begin{proposition}\label{std isotropic vectors}
Let $l\in I^{1}(L)$ be an arbitrary element.
Then there is a canonical bijection 
\[
O(L)\backslash I^{1}(L) 
\simeq  
\mathcal{G}(l^{\perp }/{\Z}l).
\]
\end{proposition}

\begin{proof}
For $l\in I^{1}(L)$,
there exists $m'\in L$ with $(l, m')=1$.
Setting $m:=m'-\frac{(m', m')}{2}l$, 
we have $(l, l)=(m, m)=0$ and  $(l, m)=1$.
Hence there is an embedding \,
$\varphi :U\hookrightarrow L$\, with\, $\varphi (f)=l$.
Since the lattice $U$ is unimodular, 
we have 
\[
L=\varphi (U)\oplus \varphi (U)^{\perp}\simeq \varphi (U)\oplus (l^{\perp }/{\Z}l) 
\] 
so that 
$A_{L}\simeq A_{l^{\perp }/{\Z}l}$ for each $l\in I^{1}(L)$.
In particular, 
$\mathcal{G}(l^{\perp }/{\Z}l)$ 
is independent of the choices of $l\in I^{1}(L)$.

We define the map 
$\mu : O(L)\backslash I^{1}(L)\to \mathcal{G}(l^{\perp }/{\Z}l)$
by 
$\mu (l'):=(l')^{\perp }/{\Z}l'$.
If there is an isometry 
$\gamma : l_{1}^{\perp }/{\Z}l_{1}\simeq l_{2}^{\perp }/{\Z}l_{2}$, 
we can extend  $\gamma $
to the isometry $\widetilde{\gamma } \in O(L)$ 
with 
$\widetilde{\gamma }(l_{1})=l_{2}$.
Therefore $l_{1}$ and $l_{2}$ are $O(L)$-equivalent  
so that 
$\mu $ is injective.
Given a lattice $K\in \mathcal{G}(l^{\perp }/{\Z}l)$, 
we have an isometry 
$\Psi :K\oplus U\stackrel{\simeq }{\rightarrow } L$
by Proposition \ref{Nikulin2}, 
which gives $\Psi (f)\in I^{1}(L)$.
Since $\mu (\Psi (f))=K$, 
$\mu $ is surjective.
\end{proof}

\begin{proposition}\label{sq}
For a primitive isotropic vector 
$l \in I^{d}(L)$, 
we have the equality 
$d^{2}\cdot \# (D_{l^{\perp}/{\Z}l})=\# (D_{L})$.
\end{proposition}

\begin{proof}
The free ${\Z}$-module 
$\widetilde{L}:=\langle L, \frac{l}{d} \rangle \subset L^{\vee }$ 
is an even overlattice of $L$. 
Since $\frac{l}{d}\in I^{1}(\widetilde{L})$, 
we have 
$\widetilde{L}\simeq U\oplus (l^{\perp}/{\Z}l)$ 
so that 
$D_{\widetilde{L}}\simeq D_{l^{\perp}/{\Z}l}$. 
If we set 
$H:=\widetilde{L}/L \subset D_{L}$, 
then $H$ is an isotropic cyclic group of order $d$ 
satisfying 
$D_{\widetilde{L}}\simeq H^{\perp}/H$. 
Hence 
$\# (D_{L})
=d^{2}\cdot \# (D_{\widetilde{L}})
=d^{2}\cdot \# (D_{l^{\perp}/{\Z}l})$. 
\end{proof}

\begin{corollary}\label{sq-free}
If ${\rm det\/}(L)$ is square-free, 
every primitive isotropic vector $l\in L$ 
satisfies ${\rm div\/}(l)=1$.
\end{corollary}

Let $l\in I^{1}(L)$. 
Choose an element $m\in L$ so that 
$(m, m)=0$ and $(l, m)=1$. 
Then we define  
$L_{0}:={\Z}l+{\Z}m \simeq U$ 
and 
$L_{1}:=L_{0}^{\perp}\simeq l^{\perp}/{\Z}l$. 
We can describe the group 
\[
O(L)^{l}:=\{ \gamma \in O(L) \, \vert \, \gamma (l)=l \}.
\]
in terms of the free Abelian group $L_{1}$ and the group $O(L_{1})$ as follows. 
For an element $v\in L_{1}$, 
define the isometry $T_{v}\in O(L)^{l}$ by 
\begin{eqnarray}
& & T_{v}(l)=l, \nonumber\\ 
& & T_{v}(m)=m+v-\frac{(v, v)}{2}l, \nonumber\\
& & T_{v}(v')=v'-(v', v)l, \: \: \: v'\in L_{1}.
 \label{eqn:translation}
\end{eqnarray}

\begin{proposition}\label{O(L)l}
For $l\in I^{1}(L)$, 
we regard $L_{1}\subset O(L)^{l}$ by 
the correspondence $v\mapsto T_{v}$ 
and regard $O(L_{1})\subset O(L)^{l}$ by 
the correspondence $g\mapsto {\rm id\/}_{L_{0}}\oplus g$.
Then  
\[
O(L)^{l}\simeq O(L_{1})\ltimes L_{1}.
\]
Moreover, 
the subgroup $L_{1}\subset O(L)^{l}$ 
acts trivially on 
$D_{L}\simeq D_{L_{1}}$.
\end{proposition}

\begin{proof}
Take an isometry $\gamma \in O(L)^{l}$. 
Since $\gamma (l^{\perp})=l^{\perp}$ 
and $l^{\perp}={\Z}l\oplus L_{1}$, 
$\gamma $ induces the isometry 
$pr_{L_{1}}\circ (\gamma |_{L_{1}}) \in O(L_{1})$. 
The correspondence 
$\gamma \mapsto pr_{L_{1}}\circ (\gamma |_{L_{1}})$ 
induces the homomorphism 
$\pi : O(L)^{l}\to O(L_{1})$. 
Then 
$\pi $ is surjective 
because 
$\pi ({\rm id\/}_{L_{0}} \oplus g)=g $ for $g\in O(L_{1})$. 
The correspondence $g\mapsto {\rm id\/}_{L_{0}} \oplus g$ 
is a section of $\pi $.

We prove that ${\rm Ker\/}(\pi )=L_{1}\subset O(L)^{l}$. 
The inclusion ${\rm Ker\/}(\pi )\supset L_{1}$ is apparent. 
For $\gamma \in {\rm Ker\/}(\pi )$, 
we can write 
$\gamma (m)=m+v+\alpha l$ 
for some integer $\alpha \in {\Z}$ and some vector $v\in L_{1}$. 
Since $(\gamma (m), \gamma (m))=(m, m)=0$, 
we can determine $\alpha $ as $\alpha =-\frac{(v, v)}{2}$. 
On the other hand, 
if we take $v'\in L_{1}$, 
we can write $\gamma (v')=v'+\beta (v')l$ for some integer $\beta (v')\in {\Z}$. 
Since $(\gamma (v'), \gamma (m))=(v', m)=0$, 
we have $\beta (v')=-(v, v')$. 
Hence we have $\gamma =T_{v}$. 
The claim that $r_{L}(L_{1})=\{ {\rm id\/}\}$ is obvious.
\end{proof}

\subsection{The Baily-Borel compactification}
In this subsection, 
we recall  
the Baily-Borel compactification of 
an arithmetic quotient of a type IV symmetric domain (\cite{B-B}, \cite{PS}),  
following \cite{Sc}. 
Although 
the Baily-Borel compactifications are defined for 
general arithmetic groups, 
we restrict ourselves to 
finite-index subgroups of $O(L)^{+}$ containing $\{ \pm {\rm id\/}\}$.

Let $L$ be an even lattice of ${\rm sign\/}(L)=(2, b^{-})$.
Set
\[ 
\Omega_{L}:=\{ {\C}\omega \in {\proj}(L_{{\C}}) \, |\, (\omega,  \omega )=0, (\omega,  \bar{\omega } )>0 \} ,
\]
which has two connected components 
$\Omega _{L}^{+}$ and $\Omega _{L}^{-}$  
exchanged by the complex conjugation on $\Omega_{L}$.
Denote by $O(L)^{+}\subset O(L)$ the subgroup of index at most $2$ 
which consists of isometries preserving $\Omega _{L}^{+}$.
By 
associating the oriented two-plane ${\R}{\rm Re\/} \omega \oplus {\R}{\rm Im\/} \omega $ 
to ${\C}\omega \in \Omega_{L}$, we have the isomorphism
\begin{eqnarray}
\Omega_{L} \simeq  
\Bigl\{ \: \Bigl. (E, \tau )\Bigr| & & E\subset L_{\R} \: \text{is a positive-definite two-plane, } \nonumber\\
& &  \tau \text{ is an orientation of } E \: \Bigr\}. 
 \label{eqn:grassmann}
\end{eqnarray}

The right hand side of (\ref{eqn:grassmann}) 
is an open subset of the oriented Grassmannian. 
Choosing a connected component of $\Omega_{L}$, 
we can attach the orientation to 
every positive-definite two-plane in $L_{{\R}}$. 
A choice of a component of $\Omega_{L}$ 
is sometimes called an {\it orientation\/} of $L$.


Regarding $\Omega _{L}^{+}$ as an open subset of 
$Q:=\{ {\C}\omega \in {\proj}(L_{{\C}})\, |\, (\omega,  \omega)=0 \} $, 
we can decompose its boundary in $Q$ as follows:
\[ 
\partial \Omega _{L}^{+}=\bigsqcup_{I} B_{I} ,
\]
where 
$I$ is an isotropic subspace of $L_{{\R}}$ and 
$B_{I}$ is the interior of 
${\proj}(I_{{\C}})\cap \overline{\Omega _{L}^{+}}$. 
Each $B_{I}$ is 
biholomorphic to the upper half-plane ${\mathbb H}$  or a point. 
The set $B_{I}$ is called a {\it rational boundary component\/} if $I=E_{{\R}}$ for some 
isotropic sublattice $E$ of $L$. 
Thus there is a canonical identification   
\[
I_{i+1}(L)\: =  
\Bigl\{ \text{ $i$-dimensional rational boundary components of } \Omega _{L}^{+} \Bigr\} 
\]
for $i=0, 1$, 
and $\Omega _{L}^{+}$ has no higher dimensional rational boundary components.

Now 
assume that 
we are given a finite-index subgroup 
$\Gamma \subset O(L)^{+}$ with $\{ \pm {\rm id\/}\}\subset \Gamma $. 
The Baily-Borel compactification of 
$\Gamma \backslash \Omega _{L}^{+} $  
is set-theoretically the set 
\begin{eqnarray*}
\overline{\Gamma \backslash \Omega _{L}^{+}}
&=& \Gamma \backslash \Bigl( \Omega _{L}^{+}\sqcup \mathop{\bigcup}_{I_{1}(L)} B_{{\Z}l} 
\sqcup \mathop{\bigcup}_{I_{2}(L)} B_{E} \Bigr) \\
&=& \Bigl( \Gamma \backslash \Omega _{L}^{+} \Bigr) 
\sqcup \mathop{\bigcup}_{\Gamma \backslash I_{1}(L)} [ {\C}l] 
\sqcup \mathop{\bigcup}_{\Gamma \backslash I_{2}(L)} \Gamma ^{E} \backslash B_{E} ,
\end{eqnarray*}
where $\Gamma ^{E}=\{ \gamma \in \Gamma \, |\, \gamma (E)=E \} $ 
for $E\in I_{2}(L)$.
It turns out that 
$\overline{\Gamma \backslash \Omega _{L}^{+}} $ is a normal projective variety. 
The boundary components of $\Gamma \backslash \Omega _{L}^{+}$
are often called {\it cusps\/}. 
We have the canonical identification 
\[
\Gamma \backslash I_{i+1}(L) = 
\Bigl\{ \text{ $i$-dimensional cusps of } \Gamma \backslash \Omega _{L}^{+} \Bigr\} .
\]

\begin{definition}
A $0$-dimensional cusp corresponding to 
$[{\Z}l]$ with ${\rm div\/}(l)=1$ 
is called a standard cusp. 
\end{definition}

By our assumption that $\Gamma \subset O(L)^{+}$, 
the notion of standardness is well-defined.  
Since 
$\{ \pm {\rm id\/}\}\subset \Gamma $, 
the set of 
$0$-dimensional standard cusps of $\Gamma \backslash \Omega _{L}^{+}$ 
is identified with 
$\Gamma \backslash I^{1}(L)$.

For a rank $2$ primitive isotropic sublattice $E\in I_{2}(L)$, 
the compact curve 
\begin{equation}
\Bigl( \Gamma ^{E} \backslash B_{E} \Bigr) \sqcup 
\Bigl( \Gamma  \backslash 
\mathop{\bigcup}_{\gamma \in \Gamma } \mathop{\bigcup}_{{\Z}l \in \gamma (E)} [{\C}l] \Bigr) 
\end{equation}
is obtained from 
the curve 
\begin{equation}
\Bigl( \Gamma ^{E} \backslash B_{E} \Bigr) \sqcup 
\Bigl( \Gamma ^{E} \backslash  \mathop{\bigcup}_{{\Z}l \in E} [{\C}l] \Bigr) \label{eqn:modular curve}
\end{equation}
by identifying $\Gamma $-equivalent points in 
$\mathop{\bigcup}_{{\Z}l \in E} [{\C}l] $.
The  curve (\ref{eqn:modular curve}) is 
the canonical compactification of 
$\Gamma ^{E} \backslash B_{E} \simeq \Gamma ^{E} \backslash {\mathbb H}$
itself.

When there exists an embedding $\varphi : U \hookrightarrow L$, 
we can write $L=\varphi (U)\oplus L_{\varphi}$ 
for the lattice $L_{\varphi}:=\varphi (U)^{\perp}\cap L$ of ${\rm sign\/}(L_{\varphi})=(1, b^{-}-1)$. 
Assume that we have chosen an orientation, 
say $\Omega _{L}^{+}$, of $L$. 
Then we can choose the connected component $L_{\varphi}^{+}$ of 
the open set 
$\{ v\in (L_{\varphi})_{{\R}} ,  (v, v)>0 \}$ 
so that 
for each vector $v\in L_{\varphi}^{+}$ 
the oriented two-plane 
${\R}\varphi (e+f) \oplus {\R}v$ belongs to  $\Omega _{L}^{+}$. 
If there is another embedding 
$\varphi ' : U \hookrightarrow L$ with $\varphi (f)= \varphi '(f)$, 
the projection 
$\gamma _{\varphi , \varphi ' }:=(pr_{L_{\varphi ' }})|_{L_{\varphi }}: L_{\varphi } \to L_{\varphi ' }$ 
is an isometry (Lemma \ref{proj}). 
Then it can be immediately checked that 
$\gamma _{\varphi , \varphi ' }$ maps 
the cone $L_{\varphi}^{+}$ to the cone $L_{\varphi '}^{+}$.

\section{${\rm FM\/}(S)$, ${\rm FM\/}^{d}(S)$ and the $0$-dimensional cusps of 
$\Gamma _{S}^{+}\backslash \Omega _{\widetilde{NS}(S)}^{+}$}

\subsection{$0$-dimensional standard cusps and FM-partners}
Let $S$ be a $K3$ surface. 
We denote by $D^{b}(S)$ 
the bounded derived category of $S$.  
A $K3$ surface $S'$ is called a 
$Fourier$-$Mukai$ $(FM)$ $partner$ of $S$ 
if $D^{b}(S)\simeq D^{b}(S')$ as triangulated categories.
Let ${\rm FM\/}(S)$ be the set of isomorphism classes of FM-partners of $S$.
By the results of Mukai and Orlov (\cite{Mu1}, \cite{Mu2}, \cite{Or}),
${\rm FM\/}(S)$ can be studied via the period of $S$. 
We denote by $\widetilde{H} (S, {\Z})$ the total cohomology group
$H^{\ast }(S, {\Z})$
equipped with the Mukai pairing 
\[
\Bigl( \, (r_{1}, l_{1}, s_{1}), \, (r_{2}, l_{2}, s_{2})\, \Bigr) :=
(l_{1}, l_{2})-(r_{1}, s_{2})-(s_{1}, r_{2})
\]
where 
$r_{i}\in H^{0}(S,{\Z}), l_{i}\in H^{2}(S,{\Z}), s_{i}\in H^{4}(S,{\Z})$.
Since 
$H^{2}(S,{\Z})\simeq \Lambda _{K3}$ 
as a lattice, 
$\widetilde{H} (S, {\Z})
\simeq H^{2}(S, {\Z})\oplus U
\simeq \widetilde{\Lambda }_{K3}$ 
as a lattice. 
We fix an isometry 
$H^{0}(S,{\Z})+H^{4}(S,{\Z})\simeq U$ 
by identifying 
$(1, 0, 0)$ with $e$, and $(0, 0, -1)$ with $f$.
The lattice $\widetilde{H} (S, {\Z})$ 
inherits the weight-two Hodge structure 
from $H^{2}(S, {\Z})$. 
If we denote by $\omega _{S}\in H^{2}(S, {\C})$ the period of $S$,  
then  
\[
\omega _{S}^{\perp}\cap \widetilde{H} (S, {\Z}) 
= \widetilde{NS}(S) 
\simeq NS(S) \oplus U.
\]
It is obvious that 
$T(S)=\widetilde{NS}(S)^{\perp}\cap \widetilde{H} (S, {\Z})$.
Let $NS(S)^{+}$ be the positive cone, 
i.e., the component of $\{ v\in NS(S)_{{\R}} ,  (v, v)>0 \}$ containing ample classes.
We choose an orientation, 
say $\Omega _{\widetilde{NS}(S)}^{+}$, 
of $\widetilde{NS}(S)$ so that 
for a vector $v\in NS(S)^{+}$ 
the oriented positive-definite two-plane ${\R}(1, 0, -1)\oplus {\R}(0, v, 0)$ 
belongs to $\Omega _{\widetilde{NS}(S)}^{+}$.

\begin{theorem}[\cite{Mu2}, \cite{Or}]\label{Mukai-Orlov}
Let $S$ and $S'$ be $K3$ surfaces. 
Then $D^{b}(S)\simeq D^{b}(S' )$ 
\, if and only if there exists a Hodge isometry 
\, $\widetilde{H} (S, {\Z})\simeq \widetilde{H} (S',  {\Z}) $.
\end{theorem}

Thanks to this theorem, 
we are able to obtain every member of ${\rm FM\/}(S)$ by an embedding of
$U$ into $\widetilde{NS}(S)$.

\begin{lemma}\label{surj}
For an embedding $\varphi \in {\rm Emb\/}(U, \widetilde{NS}(S))$,
there exists a unique (up to isomorphism) $K3$ surface
$S_{\varphi }\in {\rm FM\/}(S)$ such that 
$H^{2}(S_{\varphi }, {\Z})$ is Hodge isometric to 
$\Lambda _{\varphi } := \varphi (U)^{\perp } \cap \widetilde{H} (S, {\Z})$.
The correspondence 
$\varphi \mapsto S_{\varphi }$ gives a surjection 
\[
{\rm Emb\/}(U, \widetilde{NS}(S)) \twoheadrightarrow {\rm FM\/}(S).
\]
\end{lemma}

\begin{proof}
The lattice  
$\Lambda _{\varphi }$ is an even unimodular lattice of signature $(3, 19)$, 
hence is isometric to the $K3$ lattice $\Lambda _{K3}$.
Since $T(S)\subset \Lambda _{\varphi }$, 
$\Lambda _{\varphi }$ inherits the period from that of 
$T(S)$.
By the surjectivity of period map (\cite{To}, \cite{K3}), 
there is a $K3$ surface $S_{\varphi }$ such that 
$H^{2}(S_{\varphi }, {\Z})$ is Hodge isometric to 
$\Lambda _{\varphi } $.
By the Torelli theorem (\cite{PS-S}, \cite{K3}), 
$S_{\varphi}$ is unique up to isomorphism.
Since 
\[
\widetilde{H} (S_{\varphi},{\Z})\simeq U\oplus \Lambda _{\varphi}\simeq 
\varphi (U)\oplus \Lambda _{\varphi}= \widetilde{H} (S,{\Z})
\]
is an isometry of lattices preserving the periods, 
$S_{\varphi}\in {\rm FM\/}(S)$ by Theorem \ref{Mukai-Orlov}.

Conversely, if $S' \in {\rm FM\/}(S)$, 
there is a Hodge isometry 
$\widetilde{\Phi} :\widetilde{H} (S',{\Z})\stackrel{\simeq}{\rightarrow}\widetilde{H} (S,{\Z})$
which induces a Hodge isometry
\[
H^{2}(S', {\Z})\simeq 
\widetilde{\Phi} ( (H^{0}(S', {\Z})+ H^{4}(S', {\Z}) ) ^{\perp }\cap \widetilde{H} (S, {\Z}).
\]
It follows that 
$S'\simeq S_{\tilde{\Phi} ( (H^{0}(S', {\Z})+ H^{4}(S', {\Z}) )}$
by the Torelli theorem.
\end{proof}

By the construction of $S_{\varphi}$, 
we can identify $T(S_{\varphi })=T(S)$, 
and $NS(S_{\varphi })=\varphi (U)^{\perp}\cap \widetilde{NS}(S)$.
We shall introduce an equivalence relation on 
${\rm Emb\/}(U,\widetilde{NS}(S))$ to make the above correspondence bijective.

\begin{definition}
Set
\[
\Gamma _{S} :=r_{\widetilde{NS}(S)}^{-1} (\lambda \circ r_{T(S)}(O_{Hodge}(T(S)))) ,
\]
where 
$O_{Hodge}(T(S))$ is the group of Hodge isometries of $T(S)$, 
$r_{\widetilde{NS}(S)} :O(\widetilde{NS}(S))\rightarrow O(A_{\widetilde{NS}(S)})$  
and $r_{T(S)}:O(T(S))\rightarrow O(A_{T(S)})$ \, \, 
are the natural homomorphisms, and 
the isomorphism
$
\lambda :O(A_{T(S)})\simeq O(A_{\widetilde{NS}(S)})
$
is induced from the isometry 
$(A_{T(S)}, -q)\simeq (A_{\widetilde{NS}(S)}, q)$ 
(cf. Proposition \ref{Nikulin1}).
\end{definition}

By the identifications
$\widetilde{NS}(S)=\widetilde{NS}(S_{\varphi })$ 
and 
$T(S)=T(S_{\varphi })$,
We have  
$\Gamma _{S}=\Gamma _{S_{\varphi}}$. 
If we denote by $O_{Hodge}(\widetilde{H}(S, {\Z}))$ 
the group of Hodge isometries of $\widetilde{H}(S, {\Z})$, 
then  
$\Gamma _{S}$ is the image of the natural homomorphism
\[
O_{Hodge}(\widetilde{H}(S, {\Z}))\rightarrow O(\widetilde{NS}(S))
\]
by Proposition \ref{Nikulin1}.

In generic case, 
the period $\omega _{S}$ is not contained in any eigenspace of 
any $\varphi \in O(T(S))-\{ \pm {\rm id\/} \}$
so that $O_{Hodge}(T(S))=\{ \pm {\rm id\/} \}$. 
By the inclusions  
\[
O(\widetilde{NS}(S))_{0} \times \{ \pm {\rm id\/}\}\subseteq \Gamma _{S} \subseteq O(\widetilde{NS}(S)), 
\]
$\Gamma _{S}$ is a finite-index subgroup of $O(\widetilde{NS}(S))$. 
Then 
$\Gamma _{S}$ acts on 
${\rm Emb\/}(U,  \widetilde{NS}(S)) $ from left as 
$\gamma (\varphi ):=\gamma \circ \varphi $ 
for $\gamma \in \Gamma _{S}$ and $\varphi \in {\rm Emb\/}(U,  \widetilde{NS}(S))$.

\begin{proposition}\label{bij}
The correspondence $\varphi \mapsto S_{\varphi }$ of Lemma \ref{surj} 
induces the bijection 
\[
\Gamma _{S}\backslash {\rm Emb\/}(U, \widetilde{NS}(S))
\simeq {\rm FM\/}(S).
\] 
\end{proposition}

\begin{proof}
It suffices to show that $S_{\varphi _{1}}\simeq S_{\varphi _{2}}$ 
if and only if there exists an isometry 
$\gamma \in \Gamma _{S} $ such that 
$\varphi _{1}=\gamma \circ \varphi _{2}$.  
By the Torelli theorem, the former condition is equivalent to the existence of a
Hodge isometry $\Phi :H^{2}(S_{\varphi _{2}},  {\Z} )
\stackrel{\simeq}{\rightarrow} H^{2}(S_{\varphi _{1}},  {\Z} )$.

Given such $\Phi $,  we get a Hodge isometry $\widetilde{\Phi } :
\widetilde{H} (S_{\varphi _{2}},  {\Z}) \stackrel{\simeq}{\rightarrow} \widetilde{H} (S_{\varphi _{1}},  {\Z} )$ 
\, by adding \, $\varphi _{1}\circ \varphi _{2}^{-1}|_{\varphi _{2}(U)} :\varphi _{2}(U)
\stackrel{\simeq}{\rightarrow} \varphi _{1}(U).$
Then, 
by Proposition \ref{Nikulin1} $(2)$, 
the Hodge isometry 
$\widetilde{\Phi } |_{T(S)} :T(S)\stackrel{\simeq}{\rightarrow} T(S) $
and the isometry 
$\widetilde{\Phi } |_{\widetilde{NS}(S)} : 
NS(S_{\varphi _{2}})\oplus \varphi _{2}(U) \stackrel{\simeq}{\rightarrow} 
NS(S_{\varphi _{1}})\oplus \varphi _{1}(U) $ 
are compatible on the discriminant form. 
It follows that $\widetilde{\Phi } |_{\widetilde{NS}(S)} \in \Gamma _{S}$. 
By the construction 
$\varphi _{1}=(\widetilde{\Phi } |_{\widetilde{NS}(S)})\circ \varphi _{2}$.

Conversely, if there exists  $\gamma \in \Gamma _{S}$ such that 
$\varphi _{1}=\gamma \circ \varphi _{2}$,  
we can find a Hodge isometry 
$g : T(S)\stackrel{\simeq}{\rightarrow} T(S) $ 
which is compatible with $\gamma $ on the discriminant group.
Then, 
again by Proposition \ref{Nikulin1} $(2)$, 
$\gamma \oplus g$ extends to 
$\widetilde{\Phi } \in O_{Hodge}(\widetilde{H} (S, {\Z}))$,  
which gives a Hodge isometry between 
$H^{2}(S_{\varphi _{2}},  {\Z})=\varphi _{2}(U)^{\perp }\cap \widetilde{H} (S,  {\Z})$\, \, and 
$H^{2}(S_{\varphi _{1}},  {\Z})=\varphi _{1}(U)^{\perp }\cap \widetilde{H} (S,  {\Z})$. 
\end{proof}

Now 
we turn to the $0$-dimensional standard cusps. 
Let 
$\Omega_{\widetilde{NS}(S)}^{+}$ 
be the bounded symmetric domain associated with the lattice $\widetilde{NS}(S)$ 
(cf. Section 2.3). 
Denote by 
$\Gamma _{S}^{+}$ the arithmetic group 
$\Gamma _{S}\cap O(\widetilde{NS}(S))^{+}$.
Since the isometry 
\:  ${\rm id\/}_{NS(S)}\oplus -{\rm id\/}_{U}\in O(\widetilde{NS}(S))_{0}$ \: exchanges 
$\Omega _{\widetilde{NS}(S)}^{+}$ and $\Omega _{\widetilde{NS}(S)}^{-}$,  
we know that 
$\Gamma _{S}^{+}$ is of index $2$ in $\Gamma _{S}$. 
Recall from Section 2.3 the canonical identification 
\[
\Gamma _{S}^{+} \backslash I^{1}(\widetilde{NS}(S)) 
=
\Bigl\{ \text{ $0$-dimensional standard cusps of } 
\Gamma _{S}^{+} \backslash \Omega_{\widetilde{NS}(S)}^{+} \Bigr\} .
\]

\begin{lemma}\label{orientation}
The natural surjection 
\[
\Gamma _{S}^{+} \backslash I^{1}(\widetilde{NS}(S)) 
\longrightarrow \Gamma _{S} \backslash I^{1}(\widetilde{NS}(S))
\]
is bijective.
\end{lemma}

\begin{proof}
For $l \in I^{1}(\widetilde{NS}(S)) $, 
there exist an embedding 
$\varphi \in {\rm Emb\/}(U, \widetilde{NS}(S))$ 
such that    
$l\in \varphi (U)$.  
Then the isometry 
${\rm id\/}_{\varphi (U)^{\perp}}\oplus -{\rm id\/}_{\varphi (U)} \in O(\widetilde{NS}(S))_{0} $ 
\: preserves 
${\Z}l$\: 
and 
interchanges \: $\Omega_{\widetilde{NS}(S)}^{+}$ with $\Omega_{\widetilde{NS}(S)}^{-}$. 
Therefore 
$\Gamma _{S}^{+}\cdot l=\Gamma _{S}\cdot l$. 
\end{proof}

\begin{theorem}\label{main1}
Let $\{ e, f\} $ be 
the canonical basis of the lattice $U$.
The map 
\[
\mu _{0}:{\rm FM\/}(S)\ni [S_{\varphi }]\longrightarrow 
[\varphi (f)]\in \Gamma _{S}\backslash I^{1}(\widetilde{NS}(S)) \: .
\]
is bijective.
\end{theorem}

\begin{proof}
The map $\mu _{0}$ is well-defined by Proposition \ref{bij}. 
Since 
$l\in I^{1}(\widetilde{NS}(S)) $ 
induces an embedding 
$\varphi \in {\rm Emb\/}(U, \widetilde{NS}(S))$ 
with 
$\varphi (f)=l$, 
$\mu _{0}$ is surjective.
If $[\varphi _{1}(f)]=[\varphi _{2}(f)]$, 
then by identifying $\varphi _{i}({\Z}f)=H^{4}(S_{\varphi _{i}}, {\Z})$, 
the isometry $id_{\widetilde{H}(S, {\Z})}$ 
induces the Hodge isometry
$\widetilde{\Phi} : \widetilde{H} (S_{\varphi _{1}}, {\Z})\simeq \widetilde{H} (S_{\varphi _{2}}, {\Z}) $ 
with
$\widetilde{\Phi} (H^{4}(S_{\varphi _{1}}, {\Z}))=H^{4}(S_{\varphi _{2}}, {\Z})$. 
Now Lemma \ref{proj} implies that 
$pr_{H^{2}(S_{\varphi _{2}}, {\Z})}\circ (\widetilde{\Phi} |_{H^{2}(S_{\varphi _{1}}, {\Z})})$ 
gives a Hodge isometry between 
$H^{2}(S_{\varphi _{1}}, {\Z})$ and $H^{2}(S_{\varphi _{2}}, {\Z})$,  
since $pr_{H^{2}(S_{\varphi _{2}}, {\Z})}$ is identity on 
$T(S_{\varphi _{2}})=\widetilde{\Phi } (T(S_{\varphi _{1}}) ) $.
\end{proof}

In this way, we have 
\[
\# {\rm FM\/}(S)=\# 
\{ \text{ 0-dimensional standard cusps of } \Gamma _{S}^{+} \backslash \Omega _{\widetilde{NS}(S)}^{+} \} .
\]

\subsection{General $0$-dimensional cusps and twisted FM-partners}

In this subsection, 
using Huybrechts-Stellari's solution of C\u ald\u araru conjecture (\cite{H-S2}), 
we study the relation between 
general $0$-dimensional cusps of $ \Gamma _{S}^{+} \backslash \Omega _{\widetilde{NS}(S)}^{+} $
and 
twisted Fourier-Mukai partners of $S$.
For twisted $K3$ surfaces, see \cite{H-S1}.

Let $(S', \alpha ')$ be a twisted $K3$ surface. 
By the natural isomorphism
\[
{\rm Br\/}(S')\simeq {\rm Hom\/}(T(S'), {\Q}/{\Z}),
\]
we identify the twisting $\alpha ' \in {\rm Br\/}(S')$ 
with a surjective homomorphism 
$ 
\alpha ':T(S')\rightarrow {\Z}/d{\Z} 
$ 
where 
$d={\rm  ord\/}(\alpha ')$.
Then ${\rm  Ker\/}(\alpha ')$ is denoted by $T(S', \alpha ')$.
A twisted $K3$ surface $(S', \alpha ')$ 
is called a $twisted$ $Fourier$-$Mukai$ $partner$ of $S$
if there is an equivalence $D^{b}(S', \alpha ')\simeq D^{b}(S)$. 
By the result of Canonaco-Stellari (\cite{C-S}), 
the equivalence is of Fourier-Mukai type.
We define 
\[
{\rm FM\/}^{d}(S)
:=
\Bigl\{ \:  (S', \alpha ') : {\rm twisted\/} \: K3 \: {\rm surface\/}, 
\: \:  D^{b}(S', \alpha ')\simeq D^{b}(S), {\rm ord\/} (\alpha ')=d \: \Bigr\} / \simeq .
\]
It is obvious that 
${\rm FM\/}^{1}(S)={\rm FM\/}(S)$.

\begin{proposition}\label{construct twisted}
There exists a map 
\[
\nu _{0}: 
\Gamma _{S}^{+}\backslash I^{d}(\widetilde{NS}(S)) 
\longrightarrow {\rm FM\/}^{d}(S). 
\]
If $d=1$, 
we have $\nu _{0}=\mu _{0}^{-1}$ on 
$\Gamma _{S}^{+}\backslash I^{1}(\widetilde{NS}(S))$.
\end{proposition}

\begin{proof}
For $l\in I^{d}(\widetilde{NS}(S))$, 
we define
\[
\widetilde{M}_{l}:= 
\Bigl\langle \frac{l}{d},  \widetilde{NS}(S)  \Bigr\rangle  
 \: \subset \widetilde{NS}(S)^{\vee} , 
\]
which is an even overlattice of $\widetilde{NS}(S)$.
Via the isometry 
\[
\lambda :=\lambda _{\widetilde{H}(S, {\Z})}:
(A_{\widetilde{NS}(S)}, q)\simeq (A_{T(S)}, -q) ,
\]
we have an isotropic element $\lambda (\frac{l}{d}) \in A_{T(S)}$ 
of ${\rm order\/}=d$. 
Set 
\[
T_{l} := 
\Bigl\langle \lambda (\frac{l}{d}),  T(S)  \Bigr\rangle  
 \: \subset T(S)^{\vee} .
\]
We have a surjective homomorphism 
\[
\alpha _{l} : T_{l}\rightarrow {\Z}/d{\Z} \: \: \: 
{\rm with\/} \: \: \: 
{\rm  Ker\/} (\alpha _{l})=T(S), \: \: \:  \alpha _{l}(\lambda (\frac{l}{d}))=\bar{1} . 
\]

Since 
the isometry $\lambda $ induces the isometry 
$\bar{\lambda }:(A_{\widetilde{M}_{l}}, q)\simeq (A_{T_{l}}, -q)$, 
we have an embedding 
$
\widetilde{M}_{l}\oplus T_{l} \hookrightarrow \widetilde{\Lambda }_{K3} 
$
with both 
$\widetilde{M}_{l}$ and $T_{l}$ 
embedded primitively. 
Since 
$\frac{l}{d}  \in  I^{1}(\widetilde{M}_{l})$, 
there exists an embedding 
$\varphi : U \hookrightarrow \widetilde{M}_{l}$ 
with $\varphi (f)=\frac{l}{d}$. 
The orthogonal complement 
$\Lambda _{\varphi} := \varphi (U)^{\perp} \cap \widetilde{\Lambda }_{K3}$ 
is isometric to the $K3$ lattice $\Lambda _{K3}$  
and has the period induced from $T_{l}$. 
Moreover, 
via the orientation of  $\widetilde{M}_{l}$ (induced from $\widetilde{NS}(S)$ ),  
we can choose the positive cone $M_{\varphi}^{+}$ 
of the lattice $M_{\varphi}:=\varphi (U)^{\perp} \cap \widetilde{M}_{l}$ 
(see Section 2.3). 
By 
the surjectivity of the period map, 
there exist a $K3$ surface $S_{\varphi}$ 
and a  Hodge isometry 
$\Phi : H^{2}(S_{\varphi}, {\Z})\stackrel{\simeq}{\to} \Lambda _{\varphi} $ 
such that 
$\Phi (NS(S_{\varphi})^{+})=M_{\varphi}^{+}$. 
Pulling back the homomorphism  
$\alpha _{l} : T_{l}\rightarrow {\Z}/d{\Z} $ by $\Phi |_{T(S_{\varphi})}$, 
we obtain a surjective homomorphism  
$\alpha _{\Phi } :  T(S_{\varphi}) \to {\Z}/d{\Z}$. 
Thus we constructed a twisted $K3$ surface $(S_{\varphi}, \alpha _{\Phi })$.

If 
there are another $K3$ surface $S_{\varphi}'$ 
and a Hodge isometry
$\Phi ': H^{2}(S_{\varphi}', {\Z})\stackrel{\simeq}{\to} \Lambda _{\varphi} $ 
such that 
$\Phi '(NS(S_{\varphi}')^{+})=M_{\varphi}^{+}$,  
we have a Hodge isometry 
$\Phi '' :=  (\Phi ')^{-1}\circ \Phi :  H^{2}(S_{\varphi}, {\Z})\stackrel{\simeq}{\to} H^{2}(S_{\varphi}', {\Z})$ 
with 
$\Phi '' (NS(S_{\varphi})^{+})=NS(S_{\varphi}')^{+}$ 
and 
$(\Phi '' |_{T(S_{\varphi})})^{\ast}\alpha _{\Phi '}=\alpha _{\Phi }$. 
Composing $\Phi ''$ with an element of the Weyl group $W(S_{\varphi}')$, 
we may assume that 
$\Phi ''$ is effective 
so that 
we have $(S_{\varphi}, \alpha _{\Phi }) \simeq (S_{\varphi}', \alpha _{\Phi '})$ 
by the Torelli theorem. 
Thus 
we constructed a twisted $K3$ surface $(S_{\varphi}, \alpha _{\varphi })$ 
from an embedding $\varphi $.

If we take another embedding 
$\varphi ' : U \hookrightarrow  \widetilde{M}_{l}$ with 
$\varphi ' (f)=\frac{l}{d}$, 
it follows from Lemma \ref{proj} and Section 2.3 
that 
the projection defines an isometry 
$\Lambda _{\varphi } \stackrel{\simeq}{\to} \Lambda _{\varphi '} $ 
which is identity on $T_{l}$ 
and 
maps 
the cone $M_{\varphi}^{+}$ to 
the cone $M_{\varphi '}^{+}$. 
So 
we have a Hodge isometry 
$\Phi : H^{2}(S_{\varphi}, {\Z}) \stackrel{\simeq}{\to} H^{2}(S_{\varphi '}, {\Z}) $ 
such that 
$\Phi (NS(S_{\varphi})^{+})= NS(S_{\varphi '})^{+}$ 
and 
$(\Phi  |_{T(S_{\varphi})})^{\ast}\alpha _{\varphi '}=\alpha _{\varphi }$. 
As above, 
we obtain an isomorphism 
$(S_{\varphi}, \alpha _{\varphi }) \simeq (S_{\varphi '}, \alpha _{\varphi '})$. 
Hence we constructed a twisted $K3$ surface 
$(S_{l}, \alpha _{l})$ from 
$l\in I^{d}(\widetilde{NS}(S))$.

Next we consider the action of the group $\Gamma _{S}^{+}$.  
For an isometry $\gamma \in \Gamma _{S}^{+}$, 
set $l':=\gamma (l)$. 
The isometry $\gamma $ extends to the isometry 
$\widetilde{\gamma }:
\widetilde{M}_{l}\stackrel{\simeq }{\to} \widetilde{M}_{l'}$. 
Take an embedding 
$\varphi : U \hookrightarrow \widetilde{M}_{l}$  
with 
$\varphi (f)=\frac{l}{d}$ 
and set 
$\varphi ' := \widetilde{\gamma } \circ \varphi : U \hookrightarrow \widetilde{M}_{l'}$, 
which satisfies 
$\varphi '(f)=\frac{l'}{d}$. 
On the other hand, 
there exists a Hodge isometry $g\in O_{Hodge}(T(S))$ 
such that 
$\lambda \circ r(\gamma )=r(g)\circ \lambda $. 
Since 
$r(g)(\lambda (\frac{l}{d}))= \lambda (\frac{l'}{d})$, 
$g$ extends to the Hodge isometry 
$\widetilde{g} : T_{l} \stackrel{\simeq}{\to} T_{l'}$ 
with 
$\widetilde{g} ^{\ast}\alpha _{l'} = \alpha _{l}$. 
Because we have 
$\bar{\lambda} \circ r(\widetilde{\gamma })=r(\widetilde{g}) \circ \bar{\lambda}$, 
the isometry 
$\widetilde{\gamma } \oplus \widetilde{g}$ 
extends to the Hodge isometry  
\[
(\varphi ' \circ \varphi ^{-1}) \oplus \Phi : 
\varphi (U) \oplus \Lambda _{\varphi} 
\stackrel{\simeq}{\to} 
\varphi '(U) \oplus \Lambda _{\varphi '}  \: \: \: 
{\rm with\/} \: \: \: 
(\Phi |_{T_{l}})^{\ast}\alpha _{l'}= \alpha _{l} . 
\] 
Since $\gamma $ preserves the orientation, 
we have 
$\Phi (M_{\varphi }^{+})=M_{\varphi '}^{+}$. 
In this way 
we obtain a Hodge isometry 
$\Phi : H^{2}(S_{l}, {\Z}) \stackrel{\simeq}{\to} H^{2}(S_{l'}, {\Z})$ 
with 
$\Phi (NS(S_{l})^{+})=NS(S_{l'})^{+}$ 
and 
$(\Phi |_{T(S_{l})})^{\ast}\alpha _{l'}=\alpha _{l}$. 
As above we have 
$(S_{l}, \alpha _{l}) \simeq (S_{l'}, \alpha _{l'})$.

Finally, we see that 
$D^{b}(S_{l}, \alpha _{l})\simeq D^{b}(S)$.
Let 
$B_{l}\in H^{2}(S_{l}, {\Q})$ 
be a B-field lift of 
$\alpha _{l}\in {\rm Br\/}(S_{l})$.
Since we have a Hodge isometry 
$T(S)\simeq T(S_{l}, B_{l})\subset \widetilde{H}(S_{l}, B_{l}, {\Z})$, 
the lattice 
$\widetilde{NS}(S_{l}, B_{l}):=
T(S_{l}, B_{l})^{\perp}\cap \widetilde{H}(S_{l}, B_{l}, {\Z})$
is isogenus to 
$\widetilde{NS}(S)$, 
so 
we have an isometry 
$\widetilde{NS}(S)\simeq \widetilde{NS}(S_{l}, B_{l})$
by Proposition \ref{Nikulin2}. 
Hence there is a Hodge embedding 
\[
\gamma \oplus g : \widetilde{NS}(S)\oplus T(S) \hookrightarrow \widetilde{H}(S_{l}, B_{l}, {\Z}).
\]
Since the natural homomorphism
$r:O(\widetilde{NS}(S))\rightarrow O(A_{\widetilde{NS}(S)})$
is surjective by Proposition \ref{Nikulin2}, 
after composing $\gamma $ with 
an element of $O(\widetilde{NS}(S))$ 
(and with ${\rm id\/}_{NS(S)}\oplus -{\rm id\/}_{U} \in O(\widetilde{NS}(S))_{0}$ if necessary), 
$\gamma \oplus g$ extends to 
an orientation-preserving Hodge isometry 
$\widetilde{H}(S, {\Z})\simeq \widetilde{H}(S_{l}, B_{l}, {\Z})$.
Then our claim follows from 
Huybrechts-Stellari's theorem (\cite{H-S2}).

The proof of the assertion that 
$\nu _{0}=\mu _{0}^{-1}$ for $d=1$ can be left to the reader.
\end{proof}

\begin{definition}
Set
\[
\mathcal{FM}^{d}(S):={\rm FM\/}^{d}(S)/\sim ,
\]
where 
$(S_{1}, \alpha _{1})\sim (S_{2}, \alpha _{2})$ 
if there exists a Hodge isometry 
$g:T(S_{1})\stackrel{\simeq}{\rightarrow}T(S_{2})$
with 
$g^{\ast}\alpha _{2}=\alpha _{1}$.
Since 
$(S_{1}, \alpha _{1})\simeq (S_{2}, \alpha _{2})$ 
implies
$(S_{1}, \alpha _{1})\sim (S_{2}, \alpha _{2})$, 
the set 
$\mathcal{FM}^{d}(S)$ is well-defined.
Denote by 
$
\pi : {\rm FM\/}^{d}(S)\rightarrow \mathcal{FM}^{d}(S)
$
the quotient map.
\end{definition}

\begin{definition}
Set
\[
I^{d}(A_{\widetilde{NS}(S)}):=
\Bigl\{ 
\Bigl. 
x\in A_{\widetilde{NS}(S)}  
\: \: \Bigr| \: \: 
q_{\widetilde{NS}(S)}(x)=0 \: \in {\Q}/2{\Z}, \:   
{\rm ord\/}(x)=d \:                
\Bigr\} .
\]
Since $\widetilde{NS}(S)=NS(S)\oplus U$,  
the map 
\[
p : I^{d}(\widetilde{NS}(S))\longrightarrow I^{d}(A_{\widetilde{NS}(S)}), \: \: \: 
l\mapsto \frac{l}{d}
\]
is surjective 
by Proposition 4.1.1 of \cite{Sc}.
\end{definition}

\begin{theorem}\label{main2}
There exists a bijective map 
$
\xi _{0} : 
r(\Gamma _{S}^{+})\backslash I^{d}(A_{\widetilde{NS}(S)})
\rightarrow \mathcal{FM}^{d}(S)
$
which fits in the following commutative diagram.
$$\CD
\Gamma _{S}^{+}\backslash I^{d}(\widetilde{NS}(S)) @>\nu _{0}>> {\rm FM\/}^{d}(S) \\
@VpVV   @VV\pi V \\
r(\Gamma _{S}^{+})\backslash I^{d}(A_{\widetilde{NS}(S)}) @>>\xi _{0}> \mathcal{FM}^{d}(S).
\endCD $$
\end{theorem}

\begin{proof} 
We define $\xi _{0}$ by $\pi \circ \nu _{0} \circ p^{-1}$.
For $x=\frac{l}{d}\in I^{d}(A_{\widetilde{NS}(S)})$ 
we have an isotropic cyclic subgroup 
$\lambda (\langle x\rangle )\subset A_{T(S)}$ 
of order $d$ 
and its generator $\lambda (x)\in A_{T(S)}$, 
which define an overlattice 
$T_{x}\supset T(S)$ 
and a surjective homomorphism 
$\alpha _{x} : T_{x}\rightarrow {\Z}/d{\Z}$ 
with 
${\rm  Ker\/} (\alpha _{x}) =T(S)$, 
$\alpha _{x}(\lambda (x))=\bar{1}\in {\Z}/d{\Z}$. 
Here 
we write $\lambda $ for 
$\lambda _{\widetilde{H}(S, {\Z})} : 
(A_{\widetilde{NS}(S)}, q) \stackrel{\simeq}{\rightarrow} (A_{T(S)}, -q)$.
The isotropic subgroup 
$\langle x\rangle \subset A_{\widetilde{NS}(S)}$ 
defines an overlattice 
$\widetilde{M}_{x}=\Bigl\langle \frac{l}{d}, \widetilde{NS}(S) \Bigr\rangle  
\supset \widetilde{NS}(S)$.
Similarly as the construction of $\nu _{0}$ in Proposition \ref{construct twisted}, 
$\lambda $ induces an isometry 
$\bar{\lambda } : (A_{\widetilde{M}_{x}}, q)\simeq (A_{T_{x}}, -q)$, 
which gives an embedding 
$\widetilde{M}_{x}\oplus T_{x}\hookrightarrow \widetilde{\Lambda} _{K3}$.
Then 
the isotropic vector
$\frac{l}{d}\in I^{1}(\widetilde{M}_{x})$ 
defines the twisted $K3$ surface 
$(S_{l}, \alpha _{x})\in {\rm FM\/}^{d}(S)$.

Since 
the definitions of $T_{x}$ and $\alpha _{x}$ 
depend only on 
$x\in I^{d}(A_{\widetilde{NS}(S)})$, 
the equivalence class $[(S_{l}, \alpha _{x})]\in \mathcal{FM}^{d}(S)$ 
is independent of 
the choices of $l\in I^{d}(\widetilde{NS}(S))$ 
with $x=\frac{l}{d}$.

Take an isometry $\gamma \in \Gamma _{S}^{+}$
and set 
$x':=r(\gamma )(x)$. 
If we choose an isotropic vector $l\in I^{d}(\widetilde{NS}(S))$ 
with $\frac{l}{d}=x$, 
then $\gamma (\frac{l}{d})=x'\in A_{\widetilde{NS}(S)}$. 
Now we have 
$(S_{l}, \alpha _{l})\simeq (S_{\gamma (l)}, \alpha _{\gamma (l)})$ 
so that 
$[(S_{l}, \alpha _{l})] = [(S_{\gamma (l)}, \alpha _{\gamma (l)})]$. 
Hence the map 
$\xi _{0}$ is well-defined.

We prove the surjectivity of $\xi _{0}$.
For a twisted FM partner 
$(S', \alpha ')\in {\rm FM\/}^{d}(S)$, 
we have a Hodge isometry 
$T(S)\simeq T(S', \alpha ')$ 
and an inclusion 
$T(S', \alpha ')\subset T(S')$ 
with 
${\rm  Ker\/} (\alpha ')=T(S', \alpha ')$.
Set 
$y:=(\alpha ')^{-1}(\bar{1})\in A_{T(S)}\simeq A_{T(S', \alpha ')}$, 
where 
$\bar{1}\in {\Z}/d{\Z}$ is the generator.
Then we have an isotropic element 
$x:=\lambda ^{-1}(y)\in I^{d}(A_{\widetilde{NS}(S)})$.
%
%
%
%
By the construction of $T_{x}$ and $\alpha _{x}$, 
we have a Hodge isometry 
$(T_{x}, \alpha _{x})\simeq (T(S'), \alpha ')$ 
so that 
$[(S_{l}, \alpha _{x})]=[(S', \alpha ')] \in \mathcal{FM}^{d}(S)$,  
where 
$l \in I^{d}(\widetilde{NS}(S))$ is such that $\frac{l}{d}=x$.

Finally, 
we prove the injectivity of $\xi _{0}$.
Suppose we have 
$[(S_{l}, \alpha _{l})]=[(S_{l'}, \alpha _{l'})] \in \mathcal{FM}^{d}(S)$.
Since there exists a Hodge isometry 
$\widetilde{g}:T(S_{l})\stackrel{\simeq}{\rightarrow} T(S_{l'})$ 
with 
$\widetilde{g}^{\ast}\alpha _{l'}=\alpha _{l}$, 
$\widetilde{g}$  maps 
$T(S)={\rm Ker\/}(\alpha _{l})\subset T(S_{l})$ 
to 
$T(S)={\rm Ker\/}(\alpha _{l'})\subset T(S_{l'})$ 
isometrically.
Setting 
$g:=\widetilde{g}|_{T(S)} \in O_{Hodge}(T(S))$, 
we have 
$r(g)\Bigl( \lambda (\frac{l}{d}) \Bigr)=\lambda \Bigl( \frac{l'}{d} \Bigr) $.
Hence 
$\lambda ^{-1} \circ r(g) \circ \lambda 
\in r(\Gamma _{S}) = r(\Gamma _{S}^{+}) $
maps 
$\frac{l}{d}$ to $\frac{l'}{d}$.
\end{proof}



The set 
$r(\Gamma _{S}^{+})\backslash I^{d}(A_{\widetilde{NS}(S)})$
is relatively easy to calculate.

\begin{corollary}\label{good case}
If there is an embedding 
$U\hookrightarrow NS(S)$, 
then 
$\nu _{0} : \Gamma _{S}^{+}\backslash I^{d}(\widetilde{NS}(S)) 
\rightarrow {\rm FM\/}^{d}(S)$ 
is bijective. 
In particular, we have the formula
\[
\# {\rm FM\/}^{d}(S)=\# 
\{ \text{ 0-dimensional cusp } [l] 
\text{ of } \Gamma _{S}^{+}\backslash \Omega _{\widetilde{NS}(S)}^{+} 
\text{ with } {\rm  div\/}(l)=d \} .
\]
\end{corollary}

\begin{proof}
Since there is an embedding 
$U\oplus U \hookrightarrow \widetilde{NS}(S)$, 
$p$ is bijective by 
Proposition 4.1.3 of \cite{Sc}.
On the other hand, 
we have an embedding 
$U\oplus U \hookrightarrow \widetilde{NS}(S) \subset \widetilde{M}_{x}$ 
so that 
Proposition \ref{Nikulin2} implies that 
$\pi $ is bijective.
\end{proof}

In this way, 
we can obtain informations about 
${\rm FM\/}^{d}(S)$ 
by studying $0$-dimensional cusps of $\Gamma _{S}^{+}\backslash \Omega _{\widetilde{NS}(S)}^{+} $ 
and vice versa.

\section{${\rm FM\/}_{ell}(S)$ and certain $1$-dimensional cusps of 
$\Gamma _{S}^{+}\backslash \Omega _{\widetilde{NS}(S)}^{+}$}

Now we study the relation between 
FM-partners with elliptic fibrations 
and 
1-dimensional cusps containing $0$-dimensional standard cusps in their closures.

\subsection{FM-partners with elliptic fibrations}

\begin{definition}
Define
\[
{\rm FM\/}_{ell}(S):=\Bigl\{ \Bigl. (S',\mathcal{O} _{S'}(C'))\, \Bigr| \, S'\in {\rm FM\/}(S),  
C' \text{ is a smooth elliptic curve on } S' \Bigr\} 
/\simeq 
\]
where $(S_{1}, \mathcal{O} _{S_{1}}(C_{1}) )\simeq (S_{2}, \mathcal{O} _{S_{2}}(C_{2}) ) $ 
if there exists an isomorphism $\varphi :S_{1}\simeq S_{2} $ such that 
$\varphi ^{\ast}\mathcal{O} _{S_{2}}(C_{2})\simeq \mathcal {O} _{S_{1}}(C_{1}) $. 
In other words, $(S_{1}, \mathcal{O} _{S_{1}}(C_{1}) ) \simeq (S_{2}, \mathcal{O} _{S_{2}}(C_{2}) ) $ 
if they are isomorphic as elliptic $K3$ surfaces.
\end{definition}

For an even lattice $L$, 
let $I(L)=\cup_{d} I^{d}(L)$ be 
the set of primitive isotropic vectors in $L$.

\begin{lemma}\label{FMell}
By associating to $\mathcal{O} _{S'}(C') $ its class in $NS(S')$, 
there exists a one-to-one correspondence between
${\rm FM\/}_{ell}(S)$ and the set 
\begin{equation}
\Bigl\{ \Bigl. (S',l')\, \Bigr| \, S'\in {\rm FM\/}(S),  l'\in I(NS(S')\, ) \Bigr\} 
/\simeq . 
\label{eqn:FMell}
\end{equation} 
In  $(\ref{eqn:FMell})$, 
$(S_{1}, l_{1})\simeq (S_{2}, l_{2}) $ 
if there exists a Hodge isometry 
$\Phi :H^{2}(S_{1},  {\Z})\simeq H^{2}(S_{2}, {\Z}) $
with $\Phi (l_{1})=l_{2}$.
\end{lemma}

\begin{proof}
It is known ( Corollary 6.1 of \cite{PS-S} ) 
that 
every primitive isotropic element of $NS(S)$ 
can be transformed 
by the action of $\{ \pm {\rm id\/} \} \times W(S) $ 
to the class of a smooth elliptic curve. 
Thus, 
each element of the set (\ref{eqn:FMell}) 
is represented 
by the divisor class of a smooth elliptic curve.  
Let 
$l_{1} \in NS(S_{1})$, 
$l_{2} \in NS(S_{2})$ 
be the classes of smooth elliptic curves. 
By the Torelli theorem, 
it suffices to show that,  
if there is a Hodge isometry 
$\Phi :H^{2}(S_{1}, {\Z})\stackrel{\simeq}{\rightarrow} H^{2}(S_{2}, {\Z}) $
with $\Phi (l_{1})=l_{2} $, 
then 
there is an effective Hodge isometry 
$\Phi ': H^{2}(S_{1}, {\Z})\stackrel{\simeq}{\rightarrow} H^{2}(S_{2}, {\Z}) $
with 
$\Phi ' (l_{1})=l_{2} $. 
Let 
$\mathcal{A}(S)_{\R}$ 
be the cone in $NS(S)_{\R}$ 
generated by the ample classes of $S$ over ${\R}_{>0}$.
Since 
$l_{i}$ is contained in the boundary of  
the positive cone $NS(S_{i})^{+}$, 
we have 
$\Phi (\mathcal{A} (S_{1})_{{\R}}) \subset NS(S_{2})^{+}$. 
Then 
we have two chambers  
$\Phi (\mathcal{A} (S_{1})_{\R}) $  and  $\mathcal{A} (S_{2})_{\R} $ 
in $NS(S_{2})^{+}$, 
and  
$\Phi (l_{1})=l_{2} $ 
is contained in the closures of both chambers.
If we take 
$x_{1}\in \Phi (\mathcal{A} (S_{1})_{{\R}}) $ and  
$x_{2} \in \mathcal{A} (S_{2})_{\R} $ ,  
only finite number of walls 
$\{ W_{i} \}_{1}^{N} $ 
intersect with 
the segment $\overline{x_{1}x_{2}} $, 
and each $W_{i}$ must passes through $l_{2}$. 
Writing 
$W_{i}\cap \overline{x_{1}x_{2}} =t_{i}x_{2}+(1-t_{i})x_{1}$, 
we may assume that 
$0<t_{1}<\dots <t_{N}<1$.
If we denote by $s_{W_{i}}$ 
the reflection with respect to $W_{i}$,  
then 
$ s_{W_{N}}  \circ  \ldots  \circ  s_{W_{1}} $ 
maps 
$\Phi (\mathcal{A} (S_{1})_{\R}) $  
to  
$\mathcal{A} (S_{2})_{\R} $ 
and fixes $l_{2}$. 
Hence 
$s_{W_{N}} \circ \ldots \circ s_{W_{1}}\circ \Phi $ 
gives the desired effective Hodge isometry $\Phi '$.
\end{proof}

In what follows, we identify the two sets in Lemma \ref{FMell}.
Let ${\Z}l$ be a degenerate lattice of rank $1$, 
and consider the degenerate lattice $U\oplus {\Z}l={\Z}e+{\Z}f+{\Z}l$ of rank $3$,  
where $(e,f)=1, (e,l)=(f,l)=(l,l)=(e,e)=(f,f)=0$. 
The proof of the following theorem is 
parallel to that of Lemma \ref{surj} and Proposition \ref{bij}.

\begin{theorem}\label{FMell emb}
For an embedding 
$\varphi \in {\rm Emb\/}( U\oplus {\Z}l,  \widetilde{NS}(S) )$ 
we have
$\varphi (l)\in \varphi (U)^{\perp}\cap \widetilde{NS}(S)=NS(S_{\varphi })$  
so that we get 
$[ (S_{\varphi }, \varphi (l)) ] \in {\rm FM\/}_{ell}(S)$.
Then the correspondence 
\[
\Gamma _{S}\backslash {\rm Emb\/}( U\oplus {\Z}l,  \widetilde{NS}(S) ) \owns [\varphi ] 
 \mapsto [ (S_{\varphi },  \varphi (l)) ] \in {\rm FM\/}_{ell}(S)
\]
is a bijection.
\end{theorem}

\begin{corollary}\label{count FMell}
The following counting formula for ${\rm FM\/}_{ell}(S)$ holds:  
\[
\# \left( {\rm FM\/}_{ell}(S) \right) =\# \left( \mathop{\bigsqcup}_{[M]} \mathop{\bigsqcup}_{[k]} 
r_{T(S)}(O_{Hodge}(T(S)))\backslash O(A_{M})/r_{M}(O(M)^{k}) \right) 
\]
where $[M]$ runs over $\mathcal{G} (NS(S))$ 
and $[k]$ runs over $O(M)\backslash I(M)=O(M)\backslash I_{1}(M)$. 
\end{corollary}

\begin{proof}
We have 
\begin{eqnarray*}
&\, &\Gamma _{S}\backslash {\rm Emb\/}( U\oplus {\Z}l,  \widetilde{NS}(S) ) \\
&=&\mathop{\bigsqcup}_{[M]}  \mathop{\bigsqcup}_{[k]} 
\Gamma _{S}\backslash \Bigl\{ \Bigl. \varphi \in {\rm Emb\/}( U\oplus {\Z}l,  \widetilde{NS}(S) ) \, \Bigr| \, 
\varphi (U)^{\perp}\simeq M,  
[\varphi (l)]= [k] \Bigr\} .
\end{eqnarray*}
Since $O(M\oplus U)$ acts transitively on each set 
$ \{ \varphi \in {\rm Emb\/}( U\oplus {\Z}l,  \widetilde{NS}(S) ) \, |\,  
\varphi (U)^{\perp}\simeq M,  
[\varphi (l)]= [k] \} $
with stabilizer $O(M)^{k}$,  
the last line  can be written as 
\begin{eqnarray*}
&=&\mathop{\bigsqcup}_{[M]}  \mathop{\bigsqcup}_{[k]} 
\Gamma _{S}\backslash O(M\oplus U)/O(M)^{k} \\
&=&\mathop{\bigsqcup}_{[M]}  \mathop{\bigsqcup}_{[k]}
r_{T(S)}(O_{Hodge}(T(S)))\backslash O(A_{M})/r_{M}(O(M)^{k}), 
\end{eqnarray*} 
where we used Proposition \ref{Nikulin2} 
and the inclusion $O(M\oplus U)_{0} \subset \Gamma _{S}$ 
in the last equality.
\end{proof}


\begin{definition}
Define
\[
{\rm FM\/}_{ell,sec}(S):=
\Bigl\{ \: \Bigl. [(S',\mathcal{O} _{S'}(C'))]\in {\rm FM\/}_{ell}(S) \: \Bigr| \: 
[C'] \in NS(S'), {\rm  div\/}([C'])=1 \: \Bigr\} .
\]
\end{definition}
Note that ${\rm  div\/}([C'])=1$ if and only if the elliptic fibration associated to 
$\mathcal{O} _{S'}(C')$  has a section.
If $S$ is an elliptic $K3$ surface admitting a section, 
then ${\rm FM\/}(S)=\{ S \}$ by Corollary 2.7 of \cite{H-L-O-Y}.
Hence 
${\rm FM\/}_{ell,sec}(S)$ 
(if not empty)
is exactly the set of isomorphism classes of elliptic fibrations on $S$ admitting a section.

\begin{corollary}\label{count FMellsec}
The following equality holds:
\begin{equation*}
\# \left( {\rm FM\/}_{ell,sec}(S) \right) =
\# \left( \mathop{\bigsqcup}_{ L\in \mathcal{G} (l^{\perp }/{\Z}l)}  
r_{T(S)}(O_{Hodge}(T(S)))\backslash O(A_{L})/r_{L}(O(L)) \right) , 
\end{equation*}
where $l$ is an arbitrary standard isotropic vector in $NS(S)$, 
and 
$A_{NS(S)}$ is identified with $A_{L}$ by the isometry $NS(S)\simeq U\oplus L$. 
\end{corollary}

\begin{proof}
By Proposition \ref{std isotropic vectors}, 
$O(NS(S))\backslash I^{1}(NS(S))$ is identified with 
$\mathcal{G}(l^{\perp }/{\Z}l)$.
On the other hand, 
$r_{NS(S)}(O(NS(S))^{l})\simeq r_{L}(O(l^{\perp}/{\Z}l))$
by Proposition \ref{O(L)l}.
Now the assertion follows from Corollary \ref{count FMell}.
\end{proof}

\subsection{${\rm FM\/}_{ell}(S)$ and certain $1$-dimensional cusps}

\begin{definition}
Set 
\begin{equation}
I_{2}^{st}(\widetilde{NS}(S)) := 
\Bigl\{ \Bigl. 
E\in I_{2}(\widetilde{NS}(S)) \: \Bigr| \: 
\text{there is $e\in \widetilde{NS}(S)$ with} (e, E)={\Z} \: 
\Bigr\} . \label{eqn:I_{2}^{st}}
\end{equation}
If we take 
$f, l \in E$ so that 
$(e, f)=1$ 
and 
${\rm Ker\/}(e, \cdot )|_{E}={\Z}l$, 
then 
$E+{\Z}e={\Z}l+{\Z}f+{\Z}e \simeq U\oplus {\Z}l$. 
We may assume that 
$e$ in (\ref{eqn:I_{2}^{st}}) is taken to be isotropic. 
\end{definition}

\begin{lemma}\label{reduce}
Let $E\in I_{2}^{st}(\widetilde{NS}(S))$. 

$(1)$\, 
If $f, f'\in E$ satisfy the equalities $(e, f)=(e, f')=1$,  
then
\[ 
( S_{{\Z}e+{\Z}f}, l) \simeq 
( S_{{\Z}e+{\Z}f'}, l). 
\]

$(2)$\, 
If $e,  e'\in \widetilde{NS}(S)$ satisfy the relations $(e, E)=(e', E)={\Z}$ 
and ${\rm Ker\/}(e,\cdot )|_{E}={\rm Ker\/}(e',\cdot )|_{E}={\Z}l$,  
then
\[
( S_{{\Z}e+{\Z}f}, l) \simeq 
( S_{{\Z}e'+{\Z}f'}, l). 
\]
Here $f'\in E$ is an arbitrary element 
satisfying $(e', f')=1$.

$(3)$\, 
If $e,  e'\in \widetilde{NS}(S)$ satisfy $(e, f)=(e', f)=1$ 
for some $f\in E$, 
then
\[
( S_{{\Z}e+{\Z}f}, l) \simeq 
( S_{{\Z}e'+{\Z}f}, l'). 
\]
Here $l$ (resp. $l'$) is a generator of ${\rm Ker\/}(e,\cdot )|_{E}$ 
(resp. ${\rm Ker\/}(e',\cdot )|_{E}$).
\end{lemma}

\begin{proof}
$(1)$\, 
Since $l\in ({\Z}e+{\Z}f)^{\perp}\cap ({\Z}e+{\Z}f')^{\perp}$, 
Lemma \ref{proj} tells that the projection 
from $({\Z}e+{\Z}f)^{\perp}$ to $({\Z}e+{\Z}f')^{\perp}$ 
gives a Hodge isometry 
$H^{2}(S_{{\Z}e+{\Z}f},{\Z})\simeq H^{2}(S_{{\Z}e+{\Z}f'},{\Z})$, 
which is identity on ${\Z}l$.

$(2)$\,   
Both $\{ l, f \} $ and $\{ l, f' \} $ are basis of $E$, 
so that 
$f'=\pm f+\alpha l $ 
for some integer 
$\alpha \in {\Z}$. 
Hence 
$(e, f)=(e, \pm f')=1$, 
and
$( S_{{\Z}e+{\Z}f}, l) \simeq  ( S_{{\Z}e+{\Z}(\pm f')}, l) \simeq ( S_{{\Z}(\pm e)+{\Z}f'}, l)$ 
by $(1)$.  
On the other hand, 
$( S_{{\Z}(\pm e)+{\Z}f'}, l) \simeq  ( S_{{\Z}e'+{\Z}f'}, l) $ 
by Lemma \ref{proj}.

$(3)$\, 
As in the proof of $(2)$, 
$l'=\pm l+\beta f $ 
for some integer 
$\beta \in {\Z} $.
If we project 
$l'\in ({\Z}e'+{\Z}f)^{\perp }$ 
to 
$({\Z}e+{\Z}f)^{\perp}$, 
the image of $l'$ is given by $\pm l$. 
Now the claim follows from Lemma \ref{proj}.
\end{proof}

Now we relate 
${\rm FM\/}_{ell}(S)$ 
to 
$\Gamma _{S}^{+}\backslash I_{2}^{st}(\widetilde{NS}(S))$, 
the set of $1$-dimensional cusps of 
$\Gamma _{S}^{+} \backslash \Omega_{\widetilde{NS}(S)}^{+}$
whose closures contain 0-dimensional standard cusps, 
via Theorem \ref{FMell emb}.
Similarly as Lemma \ref{orientation}, the projection  
\[
\Gamma _{S}^{+} \backslash I_{2}^{st}(\widetilde{NS}(S)) 
\longrightarrow \Gamma _{S} \backslash I_{2}^{st}(\widetilde{NS}(S))
\]
is bijective.

\begin{definition}
Define the map 
$\mu _{1}:{\rm FM\/}_{ell}(S)\longrightarrow \Gamma _{S}\backslash I_{2}^{st}(\widetilde{NS}(S))$
by 
\[
\mu _{1}((S_{\varphi }, \varphi (l))):= [ {\Z}\varphi (f) \oplus {\Z}\varphi (l) ].
\]
\end{definition}

By the definition of $I_{2}^{st}(\widetilde{NS}(S))$, 
$\mu _{1}$ is surjective.
It follows that 
\begin{eqnarray*}
\# \Bigl( {\rm FM\/}_{ell}(S) \Bigr) 
&\geq  & \# \Bigl\{ \text{1-dimensional cusp of } \Gamma _{S}^{+} \backslash \Omega _{\widetilde{NS}(S)}^{+}
 \text{\, whose } \\
&\; & \: \: \; \text{closure contains a 0-dimensional standard cusp }  \Bigr\} .
\end{eqnarray*}

\begin{proposition}\label{inj over FMellsec}
The map
\[
\mu _{1} : 
\mu _{1}^{-1} \Bigl( \mu _{1}({\rm FM\/}_{ell,sec}(S)) \Bigr) 
\longrightarrow 
\mu _{1}({\rm FM\/}_{ell,sec}(S))
\]
is bijective. 
In particular, we have 
$\mu _{1}^{-1} \Bigl( \mu _{1}({\rm FM\/}_{ell,sec}(S)) \Bigr) ={\rm FM\/}_{ell,sec}(S)$.
\end{proposition}

\begin{proof}
Take 
$[E]=\mu _{1}((S_{{\Z}e+{\Z}f},l))
\in \mu _{1}({\rm FM\/}_{ell,sec}(S))$. 
We can find an element 
$m\in ({\Z}e+{\Z}f)^{\perp}\cap \widetilde{NS}(S)$
with $(m, m)=0, (l, m)=1$.
Write 
$U_{1}:={\Z}e+{\Z}f$ 
and 
$U_{2}:={\Z}l+{\Z}m$.
Let 
$e', f', l'\in \widetilde{NS}(S)$
be such that 
$[E]=\mu _{1}((S_{{\Z}e'+{\Z}f'},l'))$. 
We may assume that 
$E={\Z}f+{\Z}l={\Z}f'+{\Z}l'$.
Then we can write 
$f'=\alpha f+\beta l, l'=\gamma f+\delta l$ 
with 
$\alpha \delta -\beta \gamma =\pm 1$.
Since 
$(\pm (\delta e-\gamma m), f')=1$ 
and 
$(\pm (\delta e-\gamma m), l')=0$, 
we get 
$(S_{{\Z}(\pm \delta e \mp \gamma m)+{\Z}f'},l')
\simeq 
(S_{{\Z}e'+{\Z}f'},l')$
by Lemma \ref{reduce} $(2)$. 
If we set 
$e'':=\pm (\delta e-\gamma m)$ 
and 
$m':=\pm (-\beta e+\alpha m)$, 
then the correspondence 
$e\mapsto e'',\: f\mapsto f',\: l\mapsto l',\: m\mapsto m'$
gives an isometry 
$\varphi \in O(U_{1}\oplus U_{2})$.
The isometry 
$\widetilde{\varphi } := 
\varphi \oplus {\rm id\/}_{(U_{1}\oplus U_{2})^{\perp}} 
\in O(\widetilde{NS}(S)) $ 
is identity on 
$A_{\widetilde{NS}(S)}\simeq A_{(U_{1}\oplus U_{2})^{\perp}}$, 
so that 
$\widetilde{\varphi }\in \Gamma _{S}$.
It follows that 
\[
(S_{{\Z}e+{\Z}f},l)
\stackrel{\widetilde{\varphi }}{\simeq}
(S_{{\Z}e''+{\Z}f'},l')\simeq
(S_{{\Z}e'+{\Z}f'},l').
\]
\end{proof}

\begin{corollary}\label{inj if sq-free}
If ${\rm det\/}NS(S)$ is square-free,  
then every $0$-dimensional cusp of
$\Gamma _{S}^{+} \backslash \Omega _{\widetilde{NS}(S)}^{+} $ 
is standard.
Every elliptic fibration on $S'\in {\rm FM\/}(S)$ (if exists) has a section.
Therefore, $\mu _{1}$ is bijective in this case.
If ${\rm rk\/}(NS(S))\geq 3$ in addition, 
$\Gamma _{S}^{+} \backslash \Omega _{\widetilde{NS}(S)}^{+} $ 
has exactly one $0$-dimensional cusp 
corresponding to $S$.
\end{corollary}

\begin{proof}
The assertions except the last one are 
consequences of Corollary \ref{sq-free} and Proposition \ref{inj over FMellsec},
while the last assertion is Corollary 2.7 of \cite{H-L-O-Y}.
\end{proof}

However, 
the map $\mu _{1} $ is not injective in general.
By Lemma \ref{reduce} $(3)$, 
it is easily seen that 
\begin{eqnarray}
\mu _{1}^{-1}\Bigl( [E] \Bigr) 
&\simeq &  
\Gamma _{S}^{E} \backslash \Bigl\{ \: \Bigl. f\in E \: \Bigr| \: {\rm div\/}(f)=1 \: \Bigr\} 
\label{eqn:mu_{1}^{-1}}
\end{eqnarray} 
for $E\in I_{2}^{st}(\widetilde{NS}(S))$, 
where 
$\Gamma _{S}^{E}=\{ \gamma \in \Gamma _{S} \, |\,  \gamma (E)=E \} $.
The right hand side of (\ref{eqn:mu_{1}^{-1}}) 
is the set of standard cusps 
appearing in the canonical compactification of the curve 
$\Gamma _{S}^{E} \backslash B_{E} \simeq  \Gamma _{S}^{E} \backslash {\mathbb H}$ 
(cf. (\ref{eqn:modular curve})).
From the identification (\ref{eqn:mu_{1}^{-1}}), 
we observe the following two facts. 
Firstly, 
for two elliptic $K3$ surfaces 
$(S_{1}, l_{1})$ and $(S_{2}, l_{2})$ 
with 
$S_{1} \not\simeq S_{2}$,  
$\mu _{1}((S_{1}, l_{1}))=\mu _{1}((S_{2}, l_{2}))$ 
if and only if 
the two distinct standard cusps 
$\mu _{0}(S_{1})$ and $\mu _{0}(S_{2})$ 
are connected by the $1$-dimensional cusp 
$\mu _{1}((S_{1}, l_{1}))$. 
Secondly, 
fixing an elliptic $K3$ surface $(S, l)$, 
we have 
\begin{equation} 
\# \Bigl\{ \Bigl. \: 
[(S, l')]\in {\rm FM\/}_{ell}(S) 
\: \Bigr| \: 
\mu _{1}((S, l'))=\mu _{1}((S, l)) 
\: \Bigr\}  
= 
{\rm mult\/}_{\mu _{0}(S)} \Bigl( \overline{\mu _{1}((S,l))} \Bigr) , 
\end{equation}
where 
$\overline{\mu _{1}((S,l))}$ is the closure of 
the $1$-dimensional cusp $\mu _{1}((S,l))$ in $\Gamma _{S}^{+}\backslash \Omega _{\widetilde{NS}(S)}^{+}$, 
and 
${\rm mult\/}_{x}(C)$ is the multiplicity of 
the curve $C$ at the point $x$. 
Hence 
${\rm FM\/}_{ell}(S)$ 
carries informations about 
canonical compactifications of certain $1$-dimensional cusps themselves.

We remark that, 
by Proposition \ref{construct twisted} and Theorem \ref{FMell emb}, 
an elliptic fibration on an FM-partner $S'\in {\rm FM\/}(S)$ 
gives rise to a twisted FM-partner of $S$. 
If the elliptic fibration has a section, 
then $S'\simeq S$, 
and the resulting twisted FM-partner is also isomorphic to $S$. 
In particular, 
untwisted FM-partner except $S$ itself 
can never be obtained in this way.

When two elliptic $K3$ surfaces are connected by a $1$-dimensional cusp, 
there certainly exists a geometric (but not so direct) relation between them. 
Let $L\in Pic(S)$ be the line bundle associated to an elliptic fibration,
and denote by $l\in NS(S)$ the class of $L$.
Define the isotropic lattice 
$E:={\Z}l\oplus H^{4}(S,{\Z})$.
Assume we are given a hyperbolic plane 
${\Z}e'+{\Z}f'\subset \widetilde{NS}(S)$ such that
$f'\in E$ and 
$(S, l)\not\simeq (S_{{\Z}e'+{\Z}f'}, l')$.
Here $l'\in E$ is a generator of ${\rm Ker\/}(\cdot , e')|_{E}$.
For brevity, we write $S'$ instead of $S_{{\Z}e'+{\Z}f'}$.
The situation is that, 
two non-isomorphic elliptic $K3$ surfaces 
$(S, l)$ and $(S', l')$ are connected by 
the boundary curve corresponding to $E$.
Let 
$\pi _{1}:S\times S'\rightarrow S$ 
and 
$\pi _{2}:S\times S'\rightarrow S'$
be the projections.
For a coherent sheaf $\mathcal{E}$ on $S\times S'$, 
we can associate the Mukai vector 
$v_{\mathcal{E}}:=ch(\mathcal{E})\sqrt{td(S\times S')}\in H^{*}(S\times S', {\Z})$.
The following proposition gives a way to 
obtain $l\in NS(S)$ from $l'\in NS(S')$.

\begin{proposition}\label{FM transform}
$(1)$ There is a coherent sheaf $\mathcal{E}$ on $S\times S'$
such that the cohomological FM transform
\[
\Phi _{\mathcal{E}}^{H}:\widetilde{H}(S', {\Z})\rightarrow \widetilde{H}(S,{\Z}), \: \: \: 
a \mapsto \pi _{1*}(v_{\mathcal{E}}\wedge \pi _{2}^{*}a)
\]
is a Hodge isometry with 
$\Phi _{\mathcal{E}}^{H}(H^{4}(S',{\Z}))={\Z}e',\:  
\Phi _{\mathcal{E}}^{H}(H^{0}(S',{\Z}))={\Z}f' $.
By the construction of $S'$, we may assume that $\Phi _{\mathcal{E}}^{H}$
gives the identification 
$H^{2}(S', {\Z}) \simeq ({\Z}e'+{\Z}f')^{\perp} \cap \widetilde{H}(S, {\Z})$.
(In fact, $\mathcal{E}$ induces an equivalence 
$D^{b}(S')\simeq D^{b}(S)$. )

$(2)$ Let $L'\in Pic(S')$ be the line bundle representing $l'\in NS(S')$.
Then the line bundle 
\[
{\rm det\/}(R\pi _{1*}(\mathcal{E}\otimes \pi _{2}^{*}L'))\otimes {\rm det\/}(R\pi _{1*}\mathcal{E})^{-1}
\: \in Pic (S)
\] 
is a multiple of $L$.
\end{proposition}

\begin{proof}

$(1)$\: 
Firstly, 
we claim that the $H^{0}(S, {\Z})$ component of $e'$ is not $0$.
Assume $e'=(0, m, s)$. 
Since 
$f'\in E={\Z}l\oplus H^{4}(S, {\Z})$, 
we can write $f'$ degreewise as 
$f'=(0, \alpha l, t)$ for some integer $\alpha $. 
Then we have $1=(e', f')=\alpha (m, l)$.
It follows that ${\rm  div\/}(l)=1$, 
so we have 
$(S, l)\simeq (S', l')$ by Proposition \ref{inj over FMellsec}. 
This contradicts the assumption.
Since $(e', e')=0$ and $(e', f')=1$, 
then according to \cite{Mu2}, \cite{Or} (or Proposition 10.10 of \cite{Hu}), 
there exists a coherent sheaf $\mathcal{E}'$ on $S\times S'$ such that
the cohomological FM transform 
$\Phi _{\mathcal{E}'}^{H}:\widetilde{H}(S', {\Z})\rightarrow \widetilde{H}(S,{\Z})$ 
is a Hodge isometry with 
$\Phi _{\mathcal{E}'}^{H}(H^{4}(S', {\Z}))={\Z}e'$.
Replacing $\{ e', f'\} $ by $\{ -e',-f'\} $ if necessary, 
we may assume that
the positive generator $v'$ of $H^{4}(S', {\Z})$ satisfies
$\Phi _{\mathcal{E}'}^{H}(v')=e'$.
Since $(e', f')=1$, 
the $H^{0}(S', {\Z})$ component of $(\Phi _{\mathcal{E}'}^{H})^{-1}(f')$ 
is equal to $-1$.
Then, as $(\Phi _{\mathcal{E}'}^{H})^{-1}(f')$ is isotropic, 
$(\Phi _{\mathcal{E}'}^{H})^{-1}(f')=-(1, [M], \frac{[M]^{2}}{2})=-ch(M)$ 
for some line bundle $M\in Pic(S')$.
Setting $\mathcal{E}:=\mathcal{E}'\otimes \pi _{2}^{*}M$, 
we have 
$\Phi _{\mathcal{E}}^{H}(v')=\Phi _{\mathcal{E'}}^{H}(v'\wedge ch(M))
=\Phi _{\mathcal{E'}}^{H}(v')=e'$ and 
$\Phi _{\mathcal{E}}^{H}(1)=\Phi _{\mathcal{E'}}^{H}(ch(M))=-f'$.

$(2)$ \: We compute 
\begin{eqnarray*}
& & c_{1}\Bigl( {\rm det\/}(R\pi _{1*}(\mathcal{E}\otimes \pi _{2}^{*}L'))\Bigr) \\
&=& ch_{1}\Bigl( R\pi _{1*}(\mathcal{E}\otimes \pi _{2}^{*}L)\Bigr) \\
&=& \text{ the } H^{2}(S, {\Z}) \text{ component of } 
ch\Bigl( R\pi _{1*}(\mathcal{E}\otimes \pi _{2}^{*}L')\Bigr) \wedge \sqrt{td(S)} \\
&=& \text{ the } H^{2}(S, {\Z}) \text{ component of } 
\Phi _{\mathcal{E}}^{H}\Bigl( ch(L')\wedge \sqrt{td(S')} \Bigr) . 
\end{eqnarray*}
The last equation is a consequence of the Grothendieck-Riemann-Roch formula 
(cf. Corollary 5.29 of \cite{Hu}).
Decomposing $ch(L')\in \widetilde{H}(S', {\Z})$ degreewise, 
we can write $ch(L')=(1, l', 0)$. 
Since $\sqrt{td(S')} =(1, 0, 1)$,
we have $ch(L')\wedge \sqrt{td(S')} =(1, l', 1)$.
Similarly, 
\[
c_{1}\Bigl( {\rm det\/}(R\pi _{1*}\mathcal{E})\Bigr)= 
\text{ the }  H^{2}(S, {\Z}) \text{ component of } 
\Phi _{\mathcal{E}}^{H}\Bigl( (1, 0, 1) \Bigr) .
\]
Hence, 
\[
c_{1}\Bigl( {\rm det\/}(R\pi _{1*}(\mathcal{E}\otimes \pi _{2}^{*}L'))\otimes {\rm det\/}(R\pi _{1*}\mathcal{E})^{-1}\Bigr) 
=\text{ the }  H^{2}(S, {\Z}) \text{ component of } \Phi _{\mathcal{E}}^{H}(l').
\]
Since the identification 
$H^{2}(S', {\Z}) \simeq ({\Z}e'+{\Z}f')^{\perp} \cap \widetilde{H}(S, {\Z})$ 
is given by $\Phi _{\mathcal{E}}^{H}$, 
$\Phi _{\mathcal{E}}^{H}(l')$ belongs to $E={\Z}l\oplus H^{4}(S, {\Z})$.
Hence the $H^{2}(S, {\Z})$ component of $\Phi _{\mathcal{E}}^{H}(l')$
lies in ${\Z}l$.
\end{proof}

Note that the statement of Proposition \ref{FM transform} is symmetric with respect to 
$(S, l)$ and $(S', l')$.
In this way, 
we can construct $L'$ (resp. $L$) 
from $L$ (resp. $L'$) 
via a certain sheaf on $S\times S'$.
(Of course, even if a sheaf $\mathcal{E}$ on $S\times S'$ 
induces a Hodge isometry 
$\widetilde{H}(S,{\Z})\simeq \widetilde{H}(S',{\Z})$,
it is not true in general that 
the first Chern class of the line bundle
${\rm det\/}(R\pi _{2*}(\mathcal{E}\otimes \pi _{1}^{*}L))\otimes {\rm det\/}(R\pi _{2*}\mathcal{E})^{-1}$
is isotropic.)
To obtain an elliptic curve on $S'$, 
we first eliminate the fixed components of the linear system
$\Bigl| {\rm det\/}(R\pi _{2*}(\mathcal{E}\otimes \pi _{1}^{*}L))^{\pm 1}
\otimes {\rm det\/}(R\pi _{2*}\mathcal{E})^{\mp 1} \Bigr| $
and then we take a connected component of generic member of the moving part.

We give an example for which $\mu _{1}$ is not injective.

\begin{example}\label{U(r)}
Let $S$ be a $K3$ surface with $NS(S)\simeq U(r),  r>2 $, 
and denote by $r=\mathop{\prod}_{i=1}^{\tau (r)}  p_{i}^{e_{i}}$ the prime decomposition of $r$. 
Assume that $S$ is generic 
so that 
$O_{Hodge}(T(S))=\{ \pm {\rm id\/}\} $. 
Then the followings hold : 

$(1)$  $\mathcal{G} (U(r))=\{ U(r)\}  $. 

$(2)$  $\# \Bigl( {\rm FM\/}(S) \Bigr) = 2^{\tau (r)-2}\varphi (r) $, 
where $\varphi $  is the Euler function. 

$(3)$ $S$  has two non-isomorphic elliptic fibrations, whose image by  $\mu _{1} $ 
are different $1$-dimensional cusps. 

$(4)$ $\# \Bigl( {\rm FM\/}_{ell}(S) \Bigr) = 2^{\tau (r)-1}\varphi (r)$. 

$(5)$ For  $E\in I_{2}^{st}(\widetilde{NS}(S)))$ , $\# \mu _{1}^{-1}(E) =\frac{\varphi (r)}{2}$. 
 
$(6)$ $\# \Bigl( \Gamma _{S}\backslash I_{2}^{st}(\widetilde{NS}(S)) \Bigr) =2^{\tau (r)}$.
\end{example}

\begin{proof}
Write $U(r)={\Z}l+{\Z}m $ with $(l,l)=(m,m)=0,\: \: (l,m)=r $.
Then $A_{U(r)}=\langle \frac{l}{r} \rangle +\langle \frac{m}{r} \rangle  
\simeq ({\Z}/r{\Z}) \oplus ({\Z}/r{\Z}$)  with  
\[
 ({{l}\over{r}}, {{l}\over{r}} ) \equiv ( {{m}\over{r}}, {{m}\over{r}} ) \equiv 0  \: \mod \: 2{\Z}, 
\hspace{5mm} ({{l}\over{r}}, {{m}\over{r}} )\equiv {{1}\over{r}} \: \mod \: {\Z}.
\] 
We have $O(U(r))=\{ {\rm id\/}, -{\rm id\/}, \iota, -\iota \} 
\simeq ({\Z}/2{\Z}) \oplus ({\Z}/2{\Z}) $, \, 
where\, 
$\iota (l)=m,\:  \iota (m)=l $.
Therefore $O(U(r))^{l} =O(U(r))^{m} =\{ {\rm id\/} \}$.

$(1)$\: Let $L$ be an even lattice with 
${\rm sign\/}(L)=(1,1),\:  A_{L}\simeq \langle \frac{l}{r} \rangle +\langle \frac{m}{r} \rangle $.
Given a basis $\{ \widetilde{l},  \widetilde{m} \} $ of $L^{\vee}$, 
we have $ r\widetilde{l},  r\widetilde{m} \in L$.
The inclusions  
$ 
{\Z}  r\widetilde{l} + {\Z}  r\widetilde{m}   \subset L 
\subset L^{\vee} = {\Z} \widetilde{l}+ {\Z} \widetilde{m}  
$ 
imply that 
${\Z}  r\widetilde{l} +{\Z}  r\widetilde{m}   = L $. 
Thus 
the Gram matrix of $L$ with respect to the basis
$\{ r\widetilde{l},  r\widetilde{m} \} $ \: 
is divisible by $r$. 
Then 
$\vert {\rm det\/}(L(\frac{1}{r} )) \vert  = \# (A_{L})/r^{2} =1$,  
$ {\rm sign\/}(L(\frac{1}{r} ))=(1,1)$, 
and $L(\frac{1}{r} )$ is even 
so that 
$L(\frac{1}{r} ) $ must be isometric to $U$.

$(2)$\: To calculate $\# {\rm FM\/}(S)$, we shall calculate $O(A_{U(r)})$. 
\begin{claim}
$O(A_{U(r)})
\simeq 
\mathop{\prod}_{i=1}^{\tau (r)} O(A_{U(p_{i}^{e_{i}})})$.
Thus $\# O(A_{U(r)})=2^{\tau (r)}\varphi (r)$.
\end{claim}
With respect to the basis $\{ \frac{l}{r}, \frac{m}{r} \}$ 
of $({\Z}/r{\Z})\oplus ({\Z}/r{\Z})$,  
\begin{eqnarray*}
O(A_{U(r)})
&=&
\Bigl\{ \Bigl. 
\left(
\begin{smallmatrix}
a & b \\
c & d
\end{smallmatrix}
\right)
\in  GL_{2}({\Z}/r{\Z}) \: \Bigr| \:  
ad+bc\equiv 1\! \! \! \! \mod r ,\: ab\equiv cd\equiv 0\! \! \! \! \mod r 
\Bigr\} \\
&=&
\prod_{i=1}^{\tau (r)}
\Bigl\{ \Bigl. 
\left(
\begin{smallmatrix}
a & b \\
c & d
\end{smallmatrix}
\right)
\in  GL_{2}({\Z}/p_{i}^{e_{i}}{\Z}) \: \Bigr| \:  
ad+bc-1\equiv ab\equiv cd\equiv 0\! \! \! \mod p_{i}^{e_{i}} 
\Bigr\} \\
&=&
\prod_{i=1}^{\tau (r)}
O(A_{U(p_{i}^{e_{i}})}),
\end{eqnarray*}
where the second equality follows from 
the Chinese Remainder theorem.
Direct calculations show that $O(A_{U(p_{i}^{e_{i}})})=
\Bigr\{ 
\left(
\begin{smallmatrix}
a & 0 \\
0 & a^{-1}
\end{smallmatrix}
\right)
\: ; \: a\in ({\Z}/p_{i}^{e_{i}}{\Z})^{\times }
\Bigr\}
\cup
\Bigr\{ \left(
\begin{smallmatrix}
0 & b \\
b^{-1} & 0
\end{smallmatrix}
\right)
\: ; \: b\in ({\Z}/p_{i}^{e_{i}}{\Z})^{\times }
\Bigr\} $,
so that $\# O(A_{U(p_{i}^{e_{i}})})$ $=2\varphi (p_{i}^{e_{i}})$.
Since $\varphi (r) =\prod  \varphi (p_{i}^{e_{i}}) $, 
the second assertion follows.

Now $(2)$ follows from 
$(1)$, Claim 4.2, and 
the counting formula for 
$\# {\rm FM\/}(S)$ (\cite{H-L-O-Y}).

$(3)$\: $NS(S)$ has exactly four primitive isotropic elements 
$\{ \pm l, \pm m\} $.
We may assume that 
$l, m$ are classes of effective divisors.
Since $W(S)=\{ {\rm id\/}\} $, 
both of $l, m$ are classes of elliptic fibrations 
by the first part of the proof of Lemma \ref{FMell}.
The only isometry of $NS(S)$ mapping $l$ to $m$ is $\iota $,
which is not compatible with any element of 
$O_{Hodge}(T(S))=\{ \pm {\rm id\/} \}$ 
on the discriminant group.
Therefore 
$(S, l)$ and $(S, m)$ are not isomorphic.
If there exists $\gamma \in \Gamma _{S}$ which sends
$\mu _{1}((S,l))=H^{4}(S,{\Z})\oplus {\Z}l$ 
to 
$\mu _{1}((S,m))=H^{4}(S,{\Z})\oplus {\Z}m$, 
$\gamma $ maps 
$\langle \frac{l}{r} \rangle $ 
to 
$\langle \frac{m}{r} \rangle $ 
on 
$A_{U(r)\oplus U}\simeq A_{U(r)}$, 
which does not coincide with $\{ \pm {\rm id\/} \} $.

The assertion $(4)$ follows from $(2)$ and $(3)$.

$(5)$ \: By $(1)$, for each $E\in I_{2}^{st}(U(r)\oplus U)$,  we can find a basis 
$\{ l_{1}, m_{1}, e_{1}, f_{1} \} $ of $U(r)\oplus U$ such that $E={\Z}f_{1}+{\Z}l_{1}$
and 
$(l_{1}, l_{1})=(m_{1}, m_{1})=(e_{1}, e_{1})=(f_{1}, f_{1})=
(l_{1}, e_{1})=(l_{1}, f_{1})=(m_{1}, e_{1})=(m_{1}, f_{1})=0, 
(l_{1}, m_{1})=r, (e_{1}, f_{1})=1$. 
Take integers $\alpha,  \beta \in {\Z} $ 
so that 
$\beta $ and $r\alpha $ 
are coprime to each other. 
Then 
there exist $\gamma,  \delta \in {\Z} $ 
such that $\beta \delta +r\alpha \gamma =1$. 
Set 
$
f_{2}:=\alpha l_{1}+\beta f_{1} \in E, 
e_{2}:=\gamma m_{1}+\delta e_{1}, 
l_{2}:=\delta l_{1}-r\gamma f_{1} \in E, 
m_{2}:=\beta m_{1}-r\alpha e_{1}$.  
We have an isometry 
$\varphi _{(\alpha,  \beta, \gamma,  \delta )}
\in O(U\oplus U(r))^{E} $ 
which maps
$e_{1} ( {\rm resp.\/}\: f_{1}, l_{1}, m_{1} ) $ 
to 
$e_{2} ( {\rm resp.\/}\: f_{2}, l_{2}, m_{2} ) $. 
On the discriminant group $A_{U(r)\oplus U}$, 
$\varphi _{(\alpha,  \beta, \gamma,  \delta )}(\frac{m_{1}}{r} )=\beta \frac{m_{1}}{r}, 
\varphi _{(\alpha,  \beta, \gamma,  \delta )}(\frac{l_{1}}{r} )=\delta \frac{l_{1}}{r} $.
Note that, 
the image of $\delta $ in ${\Z}/r{\Z}$ is 
uniquely determined by $\alpha , \beta $ as above.
Conversely, 
given $f_{2}=\alpha l_{1}+\beta f_{1} \in E $ 
with ${\rm  div\/}(f_{2})=1$,  
we can find 
$e_{2}=\gamma m_{1}+\delta e_{1}+\epsilon f_{1}+\zeta l_{1} $ 
with 
$(f_{2}, e_{2})=1, (e_{2}, e_{2})=0 $. 
We have 
$\beta \delta +r\alpha \gamma =1$. 
Then Lemma \ref{reduce} $(3)$ assures that 
we are allowed to take $e_{2}':=\gamma m_{1}+\delta e_{1} $
instead of $e_{2} $. 
In this case, 
$l_{2}=\delta l_{1}-r\gamma f_{1}$. 
Setting 
$m_{2}:=\beta m_{1}-r\alpha e_{1}$,  
we get 
$\varphi _{(\alpha,  \beta, \gamma,  \delta )}\in O(\widetilde{NS}(S)) $ 
which maps 
$f_{1} ({\rm resp.\/}\: e_{1}, l_{1},  m_{1} )$ 
to 
$f_{2} ({\rm resp.\/}\: e_{2}',  l_{2},  m_{2} )$. 
Therefore, 
each ${\Z}e_{2}+{\Z}f_{2}+{\Z}l_{2} \in \mu _{1}^{-1}(E)$ 
can be written as 
$\varphi _{(\alpha,  \beta, \gamma,  \delta )}({\Z}e_{1}+{\Z}f_{1}+{\Z}l_{1})$
for some 
$\alpha, \beta, \gamma, \delta \in {\Z} $
with 
$\beta \delta +r\alpha \gamma =1$.
Recall that 
$\varphi _{(\alpha,  \beta, \gamma,  \delta )} 
\in \Gamma _{S}\cdot \varphi _{(\alpha ',  \beta ', \gamma ',  \delta ')}$
if and only if 
$r_{U(r)\oplus U}(\varphi _{(\alpha,  \beta, \gamma,  \delta )})=
\pm r_{U(r)\oplus U}(\varphi _{(\alpha ',  \beta ', \gamma ',  \delta ')})$.
Since 
$r_{U(r)\oplus U}(\varphi _{(\alpha,  \beta, \gamma,  \delta )})$
is expressed as the matrix 
$\left( 
\begin{smallmatrix}
\beta & 0 \\
0 & \delta
\end{smallmatrix}
\right) $, 
we see that 
$\Gamma _{S}\varphi _{(\alpha,  \beta, \gamma,  \delta )}
=\Gamma _{S}\varphi _{(\alpha',  \beta', \gamma',  \delta' )}$
if and only if
$\beta \equiv \pm \beta ', \delta \equiv \pm \delta '\mod \: r{\Z}$.
It follows that
\begin{eqnarray*}
\mu _{1}^{-1}(E) &\simeq & 
\Bigl\{ \Bigl. [\beta ], [\gamma ] \in {\Z}/r{\Z} \: \Bigr| \:  [\beta ][\gamma ] \equiv 1 \in {\Z}/r{\Z} \Bigr\}
/\{ \pm {\rm id\/} \} \\
&\simeq & ({\Z}/r{\Z})^{\times} / \{ \pm {\rm id\/} \} .
\end{eqnarray*}

The assertion $(6)$ follows from $(4)$ and $(5)$.
\end{proof}

\section{Remarks on $\Gamma _{S}^{+} \backslash \Omega _{\widetilde{NS}(S)}^{+} $}

In this section, 
we mention two remarks on the modular variety 
$\Gamma _{S}^{+} \backslash \Omega _{\widetilde{NS}(S)}^{+} $.

For $\varphi \in {\rm Emb\/}(U, \widetilde{NS}(S))$, 
the splitting $\widetilde{NS}(S)=NS(S_{\varphi })\oplus \varphi (U)$ 
induces the realization of $\Omega _{\widetilde{NS}(S)}^{+}$ as a tube domain
\[
\iota _{\varphi } : 
NS(S_{\varphi })_{{\R}}+\sqrt{-1}NS(S_{\varphi })^{+}
\stackrel{\simeq }{\longrightarrow } 
\Omega _{\widetilde{NS}(S)}^{+}, \: \: \: \:  
y\mapsto {\C}\Bigl( \varphi (e)+y-\frac{(y, y)}{2}\varphi (f)\Bigr) .
\]
Let us observe briefly the action of \, 
$O(\widetilde{NS}(S))^{\varphi (f)}\cap \Gamma _{S}^{+}$ \;
on the tube domain 
$NS(S_{\varphi })_{{\R}}+\sqrt{-1}NS(S_{\varphi })^{+}$. 
By Proposition  \ref{O(L)l} we have  
\[
O(\widetilde{NS}(S))^{\varphi (f)}\simeq 
O(NS(S_{\varphi}))\ltimes NS(S_{\varphi}),
\] 
and the action of $NS(S_{\varphi})$ is 
identity on the discriminant group.
Take an element 
$\alpha +\sqrt{-1}\omega \in NS(S_{\varphi })_{{\R}}+\sqrt{-1}NS(S_{\varphi })^{+}$. 
According to the equations (\ref{eqn:translation}), 
$m\in NS(S_{\varphi})$ 
and 
$g\in O(NS(S_{\varphi}))$ 
act as 
\begin{eqnarray*}
m(\alpha +\sqrt{-1}\omega )&=&\alpha +m+\sqrt{-1}\omega , \\
g(\alpha +\sqrt{-1}\omega )&=&g(\alpha )+\sqrt{-1}g(\omega ).
\end{eqnarray*} 
There is a normal subgroup 
$W(S_{\varphi})\ltimes NS(S_{\varphi})\subset 
O(\widetilde{NS}(S))^{\varphi (f)}\cap \Gamma _{S}^{+}$,  
where 
$W(S_{\varphi})$ is the Weyl group of $NS(S_{\varphi})$.
Since 
$\overline{\mathcal{A}(S_{\varphi})_{\R}}$ 
is the fundamental domain for the action of 
$W(S_{\varphi})$ 
on 
$NS(S_{\varphi })^{+}$, 
\begin{equation}
\Bigl( NS(S_{\varphi })_{{\R}}/NS(S_{\varphi })\Bigr) 
+\sqrt{-1}\, \overline{\mathcal{A}(S_{\varphi})_{\R}}
 \label{eqn:fundamental domain}
\end{equation}
is the fundamental domain for the action of
$W(S_{\varphi})\ltimes NS(S_{\varphi})$ 
on 
$NS(S_{\varphi })_{{\R}}+\sqrt{-1}NS(S_{\varphi })^{+}$.
Hence   
the fundamental domain for the action of 
$O(\widetilde{NS}(S))^{\varphi (f)}\cap \Gamma _{S}^{+}$ 
on 
$NS(S_{\varphi })_{{\R}}+\sqrt{-1}NS(S_{\varphi })^{+}$
is a quotient of 
(\ref{eqn:fundamental domain}).

Next, we explain 
$\Gamma _{S}^{+} \backslash \Omega _{\widetilde{NS}(S)}^{+}$ 
in connection with $Stab(S)$, 
the space of numerical locally finite stability conditions on $D^{b}(S)$, 
following Bridgeland (\cite{Br1}, \cite{Br2}).
The space $Stab(S)$ admits 
a right action of $\widetilde{GL}_{2}^{+}({\R})$, 
the universal covering of $GL_{2}^{+}({\R})$, 
and a left action of $Aut(D^{b}(S))$.
These two actions commute.
By definition, the central charge of $\sigma \in Stab(S)$ takes the form
\[
Z:E\mapsto \Bigl( \pi (\sigma ), ch(E)\sqrt{Td(S)} \Bigr) \: ,\:  E\in K(D^{b}(S))
\]
for some vector $\pi (\sigma )\in \widetilde{NS}(S)_{{\C}}$. 
The correspondence 
$\sigma \mapsto \pi (\sigma )$
defines a continuous map
$\pi :Stab(S)\rightarrow \widetilde{NS}(S)_{{\C}}$.
Set
\[
P(S):=\Bigl\{ \Bigl. \omega \in \widetilde{NS}(S)_{{\C}} \: \Bigr| \:  
{\R}{\rm Re\/} \omega +{\R}{\rm Im\/} \omega \: \text{is a positive-definite two-plane } \Bigr\} .
\]
Via the isomorphism (\ref{eqn:grassmann}), 
$P(S)$ is a principal $GL_{2}^{+}({\R})$-bundle over $\Omega _{\widetilde{NS}(S)}$ 
by the projection
\[
p : P(S)\rightarrow \Omega _{\widetilde{NS}(S)}, \: \: 
\omega \mapsto {\R}{\rm Re\/} \omega +{\R}{\rm Im\/} \omega .
\]
Here the right action of $GL_{2}^{+}({\R})$ is given by
\[
\omega \cdot g :=
(a{\rm Re\/} \omega +b{\rm Im\/} \omega )+\sqrt{-1}(c{\rm Re\/} \omega +d{\rm Im\/} \omega ), 
\] 
where 
$\omega \in P(S)$ and 
$g=
\left(
\begin{smallmatrix}
a & b \\
c & d
\end{smallmatrix}
\right)^{-1}
\in GL_{2}^{+}({\R})$.
Each tube domain realization $\iota _{\varphi }$
induces a section 
$\Omega _{\widetilde{NS}(S)}\rightarrow P(S)$.
Denote by $P^{+}(S)\subset P(S)$ the component lying over $\Omega _{\widetilde{NS}(S)}^{+}$, 
and set 
\[
\Delta (S):=\{ \delta \in \widetilde{NS}(S)\, |\, (\delta , \delta )=-2 \}.
\]
As an immediate corollary of Bridgeland's results, 
we have the following description of the modular variety
$\Gamma _{S}^{+} \backslash \Omega _{\widetilde{NS}(S)}^{+}$.

\begin{proposition}\label{Bridgeland}
There is a connected component 
$Stab^{\dag }(S)\subset Stab(S)$ 
and a subgroup $Aut^{\dag }(D^{b}(S))\subset Aut(D^{b}(S))$
such that $p\circ \pi $
induces the isomorphism
\[
Aut^{\dag }(D^{b}(S))\backslash 
Stab^{\dag }(S)/\widetilde{GL}_{2}^{+}({\R}) 
\simeq
\Gamma _{S}^{+} \backslash 
\Bigl( \Omega _{\widetilde{NS}(S)}^{+}-\bigcup_{\delta \in \Delta (S)}\delta ^{\perp} \Bigr) .
\]
\end{proposition}

\begin{proof}
According to \cite{Br1} and \cite{Br2}, 
$\pi $ induces the isomorphism
\[
Aut^{\dag }(D^{b}(S))\backslash Stab^{\dag }(S)
\simeq 
\Gamma _{S}^{+} \backslash 
\Bigl( P^{+}(S)-\bigcup_{\delta \in \Delta (S)}\delta ^{\perp} \Bigr) .
\]
Since $\pi $ is $\widetilde{GL}_{2}^{+}({\R})$-equivariant, we have
\begin{eqnarray*}
Aut^{\dag }(D^{b}(S))\backslash 
Stab^{\dag }(S)/\widetilde{GL}_{2}^{+}({\R}) 
&\simeq &
\Gamma _{S}^{+} \backslash 
\Bigl( P^{+}(S)-\bigcup_{\delta \in \Delta (S)}\delta ^{\perp} \Bigr)
/GL_{2}^{+}({\R}) \\
&\simeq &
\Gamma _{S}^{+} \backslash 
\Bigl( \Omega _{\widetilde{NS}(S)}^{+}-\bigcup_{\delta \in \Delta (S)}\delta ^{\perp} \Bigr) .
\end{eqnarray*}
\end{proof}

\section*{Acknowledgments}
The author would like to express his gratitude to
Professor K.-I. Yoshikawa for many helpful advises, discussions, and teaching.
The relation between hyperbolic planes and standard cusps was taught by him.
The author is grateful to Professor S.Hosono for valuable comments.
The author would also like to thank 
Professors T.Bridgeland and D.Kaledin and 
Doctors S.Iida and S.Ookubo 
for answering warmly to his questions.

\end{document}